\documentclass[11pt,leqno]{amsart}
\usepackage{amsmath,amsfonts,amssymb,amscd,amsthm,amsbsy,upref}
\numberwithin{equation}{section}
\newtheorem{thm}{Theorem}[section]
\newtheorem{lem}[thm]{Lemma}
\newtheorem{cor}[thm]{Corollary}
\newtheorem{prop}[thm]{Proposition}

\theoremstyle{definition}
\newtheorem*{defn}{Definition}
\newtheorem{example}[thm]{Example}
\newtheorem{prob}[thm]{Problem}
\newtheorem{quest}[thm]{Question}
\newtheorem{remark}[thm]{Remark}

\newtheorem*{rem}{Remark}
\def\EE{{\mathbb E}}
\def\nat{{\mathbb N}}
\def\R{{\mathbb R}}

\def\real{{\mathbb R}}
\def\zed{{\mathbb Z}}
\def\ds{\displaystyle}
\def\E{{\mathcal E}}
\def\X{{\mathcal X}}
\def\ep{\varepsilon}
\def\flambda{f_{(\lambda)}}
\def\dim{\operatorname{dim}}
\def\dist{\operatorname{dist}}

\def\NORM#1{{|\!|\!| #1 |\!|\!|}}
\def\chix{{\raise.5ex\hbox{$\chi$}}}
\begin{document}
\allowdisplaybreaks

\title{Systems formed by translates of one element in $L_p(\real)$}
\author[E. Odell, B. Sari,  Th. Schlumprecht \and B. Zheng]
{E. Odell$^*$, B. Sari,  Th. Schlumprecht$^*$ \and B. Zheng$^*$}
\address{Department of Mathematics\\ The University of Texas at Austin\\
1 University Station C1200\\ Austin, TX 78712-0257}
\email{odell@math.utexas.edu, btzheng@math.utexas.edu}
\address{Department of Mathematics, University of North Texas, Denton, TX 76203-5017}
\email{bunyamin@unt.edu}
\address{Department of Mathematics\\ Texas A\&M University\\
College Station, TX 77843-3368}
\email{schlump@math.tamu.edu}
\thanks{${}^*$ Research partially supported by the National Science Foundation.}
\begin{abstract}
Let $1\le p <\infty$, $f\in L_p(\real)$ and $\Lambda\subseteq \real$.
We consider the closed subspace of $L_p(\real)$, $X_p (f,\Lambda)$,
generated by the set of translations $f_{(\lambda)}$ of $f$ by $\lambda
\in\Lambda$.
If $p=1$ and $\{f_{(\lambda)} :\lambda\in\Lambda\}$ is a bounded minimal
system in $L_1(\real)$, we prove that $X_1 (f,\Lambda)$ embeds almost
isometrically into $\ell_1$.
If $\{f_{(\lambda)} :\lambda\in\Lambda\}$ is an unconditional basic
sequence in $L_p(\real)$, then $\{f_{(\lambda)} : \lambda\in\Lambda\}$
is equivalent to the unit vector basis of $\ell_p$ for $1\le p\le 2$ and
$X_p (f,\Lambda)$ embeds into $\ell_p$ if $2<p\le 4$.
If $p>4$, there exists $f\in L_p(\real)$ and $\Lambda \subseteq \zed$
so that $\{f_{(\lambda)} :\lambda\in\Lambda\}$ is unconditional basic
and $L_p(\real)$ embeds isomorphically into $X_p (f,\Lambda)$.
\end{abstract}
\maketitle

\baselineskip=18pt      

\section{Introduction}

Let $f:\real\to\real$ and $\lambda\in\real$.
We denote by
$\flambda$ the translation of $f$ $\lambda$-units to the right for $\lambda>0$
(and $|\lambda|$-units to the left for $\lambda <0$).
Precisely,
$$\flambda (x) = f(x-\lambda)\ \text{ for }\ x\in \real\ .$$
If $f\in L_p (\real)$, $1\le p<\infty$ and $\Lambda \subseteq \real$, we
let $X_p(f,\Lambda)$ equal $[\{\flambda :\lambda\in\Lambda\}]$, where
$[\cdot]$ denotes the closed linear span in $L_p(\real)$.
Our main focus shall be on the nature of such subspaces given that
$\{\flambda :\lambda\in\Lambda\}$ has some additional structure and
$\Lambda$ is {\em uniformly discrete\/}, i.e.,
$$\inf \{ |\lambda-\lambda'| : \lambda,\lambda' \in \Lambda,\
\lambda\ne\lambda'\} >0\ .$$
The ``additional structure'' takes several forms:
$\{ \flambda :\lambda\in\Lambda\}$ is a bounded minimal system, or is
unconditional basic or can be ordered to be a (Schauder) basis or a
(Schauder) frame for $X_p (f,\Lambda)$.
It is worth mentioning that it is known that if
$\{\flambda :\lambda\in\Lambda\}$ is a bounded minimal system, in particular,
if it can be ordered to be basic, then $\Lambda$ must be uniformly discrete.
This is easy (Proposition~\ref{prop1.1} below).

The nature of $X_p (f,\Lambda)$ and $\{\flambda :\lambda\in\Lambda\}$
have been studied in a number of papers, mainly using techniques of
harmonic analysis.
Our techniques will be, largely, from the geometry of Banach spaces.
We recall seven theorems, beginning with Wiener's famous Tauberian theorem.

\begin{thm}\label{thm:Wi}
\cite{Wi}.
For $f\in L_2(\real)$, $X_2 (f,\real) = L_2(\real)$ if and only if
$\hat f(t) \ne0$ a.e.
For $f\in L_1 (\real)$, $X_1 (f,\real) = L_1(\real)$ if and only if
$\hat f(t)\ne0$ for all $t\in\real$.
\end{thm}

\begin{thm}\label{thm:AO}
\cite[Theorem 2.1]{AO}.
Let $2<p<\infty$.
There exists $f\in L_p(\real)$, all of whose derivatives exist and are in
$L_2(\real)$ (i.e., $f\in H^{2,\infty}(\real)$) so that $X_p (f,\zed) =
L_p(\real)$.
Moreover, $f$ can be chosen to satisfy, in addition, any one of the
following conditions.
\begin{itemize}
\item[(1)] $X_p (f,\nat_0) = L_p (\real)$.
\item[(2)] $(f_{(n)})_{n\in\zed}$ is orthogonal in $L_2(\real)$.
\item[(3)] $(f_{(n)})_{n\in\zed}$ is a bounded minimal system.
\end{itemize}
\end{thm}

\begin{thm}\label{thm:AO-b}
\cite{AO}.
Let $1\le p\le 2$ and let $F\subseteq L_p(\real)$ be a finite set.
Then $[\{ f_{(n)} : f\in F$, $n\in\zed\}] \ne L_p(\real)$.
\end{thm}

\begin{thm}\label{thm:ER}
\cite[Corollary 2.11]{ER}.
Let $1\le p<\infty$, $0\ne f\in L_p(\real)$.
Then $\{f_{(\lambda)} : \lambda\in \real\}$ is linearly independent.
\end{thm}

\begin{thm}\label{thm:Ol}
\cite{Ol}.
Let $\Lambda = \{\lambda_n\}_{n\in\zed} \subseteq \real$ so that
$\Lambda\cap \zed=\emptyset$ and $\lim_{|n|\to\infty} |\lambda_n - n|=0$.
Then there exists $f\in L_2 (\real)$ so that $X_2 (f,\Lambda) = L_2(\real)$.
\end{thm}

\begin{thm}\label{thm:OZ}
\cite[Theorem 2]{OZ}.
There is no unconditional basis of translates of $f$, $\{\flambda :\lambda
\in\Lambda\}$, with $X_2 (f,\Lambda) = L_2(\real)$.
\end{thm}

So the space $X_p (f,\Lambda)$, $\Lambda$ uniformly discrete, can equal
$L_p(\real)$, for $p\ge 2$ at least.
For $p=1$ the situation is different as pointed out to us by J.~Bruna.

\begin{thm}\label{thm:folk}
Let $f\in L_1(\real)$ and let $\Lambda\subseteq \real$ be uniformly discrete.
Then $X(f,\Lambda) \ne L_1(\real)$.
\end{thm}

This seems to be a folk theorem and we were unable to find a
reference. It follows from Theorem~\ref{thm:Wi} and
Lemma~\ref{lem:minimal} below (and can also be deduced from
\cite{BOU} and the proof of Lemma~\ref{lem:minimal}).

For $1<p<\infty$ it remains an open problem whether there exists
$\Lambda\subseteq \real$ and $f\in L_p(\real)$ so that, in some order,
$\{\flambda :\lambda\in\Lambda\}$ is a basis for $L_p (\real)$.

In section 2 we prove that if $f\in L_1(\real)$ and $\{\flambda :\lambda
\in \Lambda\}$ is a bounded minimal system for $X_1 (f,\Lambda)$, then
$X_1 (f,\Lambda)$ embeds almost isometrically into $\ell_1$.
The same conclusion holds if $\Lambda$ is uniformly discrete and
$\{\flambda :\lambda\in\Lambda\}$ can be ordered to be a (Schauder) frame
for $X_1(f,\Lambda)$.

In Corollary~\ref{cor2.10} we show that for $1\le p\le 2$, if
$(\flambda)_{\lambda\in\Lambda}$ is an unconditional basic sequence
then $(\flambda)_{\lambda\in\Lambda}$ is equivalent to the unit vector
basis of $\ell_p$.
For $2<p\le 4$ we show (Theorem~\ref{thm2.11})
that if $(\flambda)_{\lambda\in\Lambda}$ is unconditional basic then
$X_p(f,\Lambda)$ embeds isomorphically into $\ell_p$
but (Proposition~\ref{prop2.14})    
$(\flambda)_{\lambda\in\Lambda}$ need
not be equivalent to the unit vector basis of $\ell_p$.
For $4<p<\infty$ we give an example (Theorem~\ref{thm2.13})
of an unconditional basic sequence $(\flambda)_{\lambda\in\Lambda}$ so that
$X_p(f,\Lambda)$ contains an isomorph of $L_p[0,1]$
(which, of course, is isometric to $L_p(\real)$).

Among further results in section~2, we also consider the translation
problem for the translation invariant
space $L_p(\real) \cap L_2 (\real)$, $2<p<\infty$, and
show (Proposition~\ref{prop:2.17bis}) that if $(\flambda)_{\lambda\in\Lambda}$
is unconditional basic then it is equivalent to the unit vector basis of
$\ell_2$.


In the beginning of section 3, we revisit the problem for integer
translates of $f\in L_p(\real)$.
We also prove that if $f\in L_1 (\real)$ with $\hat f(t)\ne0$ for all $t$,
then $X_1 (f,\zed)$ embeds into $\ell_1$ (Proposition~\ref{prop3.4bis}).
We also consider discrete versions of our problem for $\ell_p (\zed,X)$
in Propositions~\ref{prop3.3}, \ref{prop3.6} and Corollary~\ref{cor3.10}.
Fourier analysis plays a role in some of these results.

In section 4, we recall some additional known results from the literature
and list some remaining open problems.

We use standard Banach space notation as may be found in \cite{LT}
or \cite{JL}.
Background material on bases, unconditional bases and such can be found there.
For the benefit of those less familiar with these notions we recall some
definitions and facts.
A {\em biorthogonal system\/} is a sequence
$(x_i,x_i^*)_{i=1}^\infty\subseteq X\times
X^*$ where $x_i^* (x_j) = \delta_{(i,j)}$.
A biorthogonal system $(x_i,x_i^*)_{i=1}^\infty\subseteq X \times X^*$
is {\em fundamental\/} (or {\em complete\/})
if $[(x_i)_{i\in\nat}] =X$ and {\em bounded\/}
if $\sup_i \|x_i\|\, \|x_i^*\| <\infty$.

$(x_i)_{i=1}^\infty \subseteq X$ is a {\em minimal system} if there exists
$(x_i^*)_{i=1}^\infty \subseteq X^*$ so that $(x_i,x_i^*)_{i=1}^\infty$
is a biorthogonal system.
This is equivalent to $x_i \notin [x_j :j\ne i]$ for all $i\in \nat$.
$(x_i)_{i=1}^\infty$ is a {\em bounded minimal system} if, in addition,
$(x_i,x_i^*)_{i=1}^\infty$ is a bounded biorthogonal system.
This is equivalent to $\inf_i \dist (x_i,[x_j :j\ne i]) >0$.
$(x_i)_{i=1}^\infty\subseteq  X$ is a (Schauder) {\em basis\/} for $X$ if
for all $x\in X$ there exists a unique sequence of scalars $(a_i)_{i=1}^\infty$
so that $x= \sum_{i=1}^\infty a_i x_i$.
This is equivalent to saying that all $x_i\ne0$, $[(x_i)]=X$ and for some
$K<\infty$, all $m<n$ in $\nat$ and all $(a_i)_1^n\subseteq \real$,
$\|\sum_{i=1}^m a_i x_i\| \le K\|\sum_{i=1}^n a_i x_i\|$.
The smallest such $K$ is the {\em basis constant\/} of $(x_i)$.
A basis $(x_i)_{i=1}^\infty$ for $X$ is
a fundamental bounded minimal system for $X$.
In this case every $x\in X$ can be written uniquely as
$x= \sum_{i=1}^\infty x_i^* (x) x_i$.
The $x_i^*$'s are a {\em basic sequence\/} in $X^*$, i.e., form a basis
for $[(x_i^*)]\subseteq X^*$ and are a basis for $X^*$ if $X$ is reflexive.
$(x_i)_{i=1}^\infty$ is an {\em unconditional basis}
 for $X$ if for all $x\in X$
there exists a unique sequence of scalars $(a_i)_{i=1}^\infty$ so that
$x= \sum_{i=1}^\infty a_i x_i$ and the {\em convergence is unconditional\/},
i.e., $x= \sum_{i=1}^\infty a_{\pi (i)} x_{\pi(i)}$ for all permutations
$\pi$ of $\nat$.
This is equivalent to all $x_i$'s $\ne 0$, $[(x_i)_{i\in\nat}] = X$ and
$$\sup \bigg\{\Big\| \sum_{i=1}^\infty \ep_i a_i x_i\Big\|
:  \sum_{i=1}^\infty a_i x_i\in B_X\ \text{ and }\
\ep_i = \pm 1\ \text{ for all }\ i\bigg\} <\infty\ .$$
Here $B_X$ denotes the closed unit ball of $X$.
This number is called the {\em unconditional basis\/} constant of
$(x_i)_{i=1}^\infty$.
The biorthogonal functionals then form an unconditional basic sequence in $X^*$.

A {\em block basis\/} $(y_i)_{i=1}^\infty$ of a basic sequence
$(x_i)_{i=1}^\infty$ is a non-zero sequence given by
$$y_i=\sum_{j=n_{i-1}+1}^{n_i} a_j x_j\ \text{ for some sequence }\
n_0 < n_1 < n_2 < \cdots$$
in $\nat_0$ and scalars $(a_j)_{j=1}^\infty \subseteq \real$.
A block basis is a basic sequence, which is unconditional basic if the $x_i$'s
are unconditional basic.
A sequence $(x_i)$ is semi-normalized if $0< \inf \|x_i\| \le \sup_i
\|x_i\| <\infty$.

A {\em Schauder frame}
for a Banach space $X$ is a sequence
$(x_i,f_i) \subseteq X\times X^*$ such that
for all $x\in X$, $x= \sum_{i=1}^\infty f_i(x) x_i$.
Of course every basis for $X$ is a frame for $X$ and just as in the basis
case, the uniform boundedness principle yields
$\sup \{ \|\sum_1^n f_i (x) x_i\| : n\in\nat,\ x\in S_X\}<\infty$
(called the {\em frame constant\/}) where
$S_X = \{x\in X: \|x\|=1\}$ is the unit sphere of $X$.
More on  frames can be found in \cite{CHL} and \cite{CDOSZ}.
Schauder frames should not be confused with Hilbert frames which are
much more restrictive.
Note that for frames, $(x_i,f_i)$ is not assumed to be a biorthogonal
sequence.

In our situation, where we are concerned with $(f_i)_{i=1}^\infty$ being
a sequence of uniformly discrete translations of some $f\in L_p(\real)$,
we do not know of an example where $(f_i)$ is a frame but is not basic.
However, many of our results would hold only given the property of
Proposition~\ref{prop2.1} below and so we have stated them in terms of frames.

Some background material on $L_p$ spaces which we shall use can be found
in \cite{AOd} and in the basic concepts chapter of \cite{JL}.
In particular we shall use that a normalized unconditional basic sequence
$(f_i)$ in $L_p(\real)$ satisfies  for constants $A_p$ and $B_p$, depending
on $p$ and the unconditional basis constant of $(f_i)$,
\begin{itemize}
\item[(1.1)] For $1\le p\le 2$, for all $(a_i) \subseteq \real$,
$$A_p^{-1} \bigg( \sum_{i=1}^\infty a_i^2\bigg)^{1/2}
\le \Big\| \sum_{i=1}^\infty a_i f_i\Big\|_p
\le B_p \bigg(\sum_{i=1}^\infty |a_i|^p\bigg)^{1/p}\ . $$
\item[(1.2)] For $2<p<\infty$, $(a_i) \subseteq \real$,
$$A_p^{-1} \bigg( \sum_{i=1}^\infty |a_i|^p\bigg)^{1/p}
\le \Big\| \sum_{i=1}^\infty a_i f_i\Big\|_p
\le B_p \bigg(\sum_{i=1}^\infty |a_i|^2\bigg)^{1/2}\ .$$
\end{itemize}

If $(f_i)$ is unconditional basic in $L_p[0,1]$, $1\le p<\infty$ then
for some $C_p$, depending on $p$ and the unconditional basis constant of
$(f_i)$, for all $(a_i)\subseteq \real$,
\begin{itemize}
\item[(1.3)] (Square function inequality)
$$\Big\| \sum_{i=1}^\infty a_i f_i\Big\|_p
\ \buildrel {C_p}\over \sim\
\left( \int_0^1 \bigg( \sum_{i=1}^\infty |a_i|^2 |f_i (t)|^2\bigg)^{p/2}
\, dt\right)^{1/p}\ .$$
\end{itemize}
Here we use ``$A\ \buildrel C\over\sim\  B$'' to denote $A\le CB$ and
$B\le CA$.

The {\em Haar basis\/} $(h_n)_{n=1}^\infty$ is a basis for $L_1[0,1]$.
This sequence is given by
\begin{equation*}
\begin{split}
(h_n)_{n=1}^\infty &= (\X_{[0,1]},\ \X_{[0,\,1/2]} - \X_{[1/2,\,1]},\
\X_{[0,\,1/4]} - \X_{[1/4,\, 1/2]}, \\
&\qquad \qquad \X_{[1/2,\,3/4]} - \X_{[3/4,\, 1]}, \
\X_{[0,\,1/8]} - \X_{[1/8,\, 1/4]}, \ldots)
\end{split}
\end{equation*}
The same system is an unconditional basis for $L_p[0,1]$, $1<p<\infty$.
Usually below, we will let $(h_n)_{n=1}^\infty$ refer to
the normalized Haar basis, i.e., $(h_n/\|h_n\|_p)_{n=1}^\infty$.
We can get an unconditional basis for $L_p(\real)$ from this by copying
$(h_n)_{n=1}^\infty$ onto each interval $[k,k+1]$, $k\in\zed$.
In this case we will have functions $(h_{n,k})_{n\in \nat,\ k\in\zed}$
and we will presume they are {\em linearly ordered\/} so as to be
{\em compatible\/} with the Haar basis ordering on each $[k,k+1]$, i.e.,
if the functions are ordered as $(x_i)_{i=1}^\infty$ and if $x_i = h_{n,k}$,
$x_j = h_{m,k}$ with $i<j$, then $n<m$.
This ordering yields that if $(g_i)_{i=1}^\infty$ is a block basis of the
Haar basis, then $(g_i|_{[n,m]})_{i=1}^\infty$ is also a block basis of
the Haar basis (well, some $g_i$'s could be 0 here) for all
integers $n<m$.

The {\em Rademacher sequence\/} $(r_n)_{n=1}^\infty$ is given by
$(r_n)_{n=1}^\infty = (h_1,h_2,h_3 + h_4, h_5 +\cdots + h_8,\ldots)$,
where the $h_n$'s refer to the non-normalized Haar functions.
It is {\em equivalent to the unit vector basis of $\ell_2$} in all $L_p[0,1]$
spaces, $1\le p<\infty$, i.e.,
$$\Big\| \sum_{i=1}^\infty a_i r_i\Big\|_p
\ \buildrel {K_p}\over\sim\
\bigg( \sum_{i=1}^\infty |a_i|^2\bigg)^{1/2}\ .$$

One reason for taking $\Lambda$ to be uniformly discrete in our considerations
is, as mentioned above, given by the easy

\begin{prop} \label{prop1.1}
\cite[Theorem 1]{OZ}.
Let $1\le p<\infty$ and let $(f_i,g_i)_{i=1}^\infty$ be a bounded
biorthogonal system in $L_p(\real)$ so that for some $f\in L_p(\real)$ and
$(\lambda_i)_1^\infty \subseteq \real$, $f_i = f_{(\lambda_i)}$ for all $i$.
Then $\Lambda = (\lambda_i)_1^\infty$ is uniformly discrete.
\end{prop}

\begin{proof}
If not, there exist subsequences $(i_m)$ and $(j_m)$ of $\nat$ so that
$\lim_{m\to\infty} |\lambda_{i_m} - \lambda_{j_m} | =0$ and
$\lambda_{i_m} \ne \lambda_{j_m}$ for all $m$.
Then
$$\|g_{i_m} \| \ge
\frac{\langle g_{i_m}, f_{i_m} - f_{j_m}\rangle} {\|f_{i_m} - f_{j_m}\|_p}
= \frac1{\|f_{i_m} - f_{j_m}\|_p}$$
and the latter is unbounded in $m$, a contradiction.
\end{proof}

\section{Main Results}\label{sec:main}

We begin with the elementary but very useful

\begin{prop}\label{prop2.1}
Let $\Lambda \subseteq \real$ be uniformly discrete, $1\le p<\infty$, and
$f\in L_p (\real)$.
Then for all intervals $I = [a,b]$,
$\sum_{\lambda\in\Lambda}\|\flambda |_I\|_p^p <\infty$.
\end{prop}

\begin{proof}
Choose $\ep_0 >0$ so that $|\lambda -\lambda'| >\ep_0$ for all
$\lambda,\lambda'\in\Lambda$ with $\lambda \ne \lambda'$.
For $\ell\in\zed$, set $I_\ell = [a+ (\ell-1) \ep_0, a+\ell \ep_0]$.
Then for $\ell\in\zed$,
\begin{equation*}
\sum_{\lambda\in\Lambda} \|\flambda |_{I_\ell}\|_p^p
 = \sum_{\lambda\in\Lambda} \int_{I_\ell} |f(t-\lambda)|^p\,dt
=\sum_{\lambda\in\Lambda} \int_{a+(\ell-1)\ep_0-\lambda}^{a+\ell\ep_0-\lambda}
|f(t)|^p\,dt \le \|f\|_p^p\ ,
\end{equation*}
since the intervals of integration are disjoint for $\lambda\in\Lambda$.
Thus
\begin{equation*}
\sum_{\lambda\in\Lambda} \|\flambda |_I \|_p^p
\le \sum_{\ell=1}^{\left\lceil \frac{b-a}{\ep_0}\right\rceil}
\sum_{\lambda\in\Lambda} \|\flambda |_{I_\ell} \|_p^p
\le \left\lceil \frac{b-a}{\ep_0}\right\rceil \|f\|_p^p\ .\qquad \qed
\end{equation*}
\renewcommand{\qed}{}
\end{proof}

We note a simple consequence of Proposition~\ref{prop2.1}.
We remark that in \cite[Theorem~4.1]{AO}, it is proved that if $1<p<\infty$
and $f\in L_p(\real) \cap L_1 (\real)$ then $X_p (f,\zed)\ne L_p(\real)$.

\begin{prop}\label{prop2.3bis}
Let $1 < p< \infty$, $f\in L_p(\real)\cap L_1(\real)$, and let
$(f_i)_{i=1}^\infty$ be a sequence of uniformly discrete translates of $f$.
Then $(f_i)_{i=1}^\infty$ is not a fundamental bounded minimal system for
$L_p (\real)$.
Furthermore, there is no sequence $(g_i)_{i=1}^\infty \subseteq L_q(\real)$
$(1/p + 1/q=1)$ so that $(f_i,g_i)_{i=1}^\infty$ is a frame for
$L_p(\real)$.
\end{prop}

\begin{proof}
Assume $(f_i,g_i)$ were in fact such a frame.
$\|f_i\|_p = \|f\|_p$ for all $i$ and thus $(g_i)_{i=1}^\infty$ is
$\omega^*$-null and hence bounded in $L_q(\real)$.
Let $K = \sup_i \|g_i\|_q$.
Choose $n_0\in\nat$ with
$$\sum_{j=n_0+1}^\infty \|f_i |_{[0,1]} \|_1 < \frac1{4K} \ .$$
Choose $h:\real\to\real$ so that $|h| = \chi_{[0,1]}$ and
$|\langle h,g_i\rangle| < \frac1{4n_0\|f\|_1}$ for $i\le n_0$ ($h$ could
be a Rademacher function).
Thus $\|h\|_p = \|h\|_1 =1$.
Also $h = \sum_{i=1}^\infty \langle h,g_i\rangle f_i$, the series
converging in $L_p(\real)$, and so
$$h|_{[0,1]} = \sum_{i=1}^\infty \langle h,g_i\rangle f_i|_{[0,1]}\ ,$$
the series converging in $L_1 [0,1]$.
Then
\begin{equation*}
\begin{split}
1& = \|h\|_1 \le \sum_{i=1}^\infty |\langle h,g_i\rangle| \|f_i|_{[0,1]}\|_1\\
\noalign{\vskip6pt}
& \le \sum_{i=1}^{n_0} |\langle h,g_i\rangle| \|f\|_1
+ \sum_{i=n_0+1}^\infty \|g_i\|_q \|f_i|_{[0,1]} \|_1\\
\noalign{\vskip6pt}
& < \frac{n_0}{4n_0 \|f\|_1} \|f\|_1
+ \sup_i \|g_i\|_q \frac1{4K} = \frac12\ ,
\end{split}
\end{equation*}
a contradiction.


The argument is similar if we assume that $(f_i,g_i)_{i=1}^\infty$ is a
fundamental bounded biorthogonal system for $L_p(\real)$.
Then, for the same $h$, $n_0$ and for $\ep>0$ arbitrary, we can
choose $f=\sum_{i=1}^n a_i f_i$ with $\|h-f\|_p < \ep$.
Thus $\|f|_{[0,1]} -h\|_1 < \ep$ and
\begin{equation}\label{eq:fbbio}
1-\ep \le \|f|_{[0,1]}\|_1
\le \sum_{i=1}^{n_0} |a_i|\, \|f_i|_{[0,1]}\|_1
+ \sum_{i=n_0+1}^n |a_i|\, \|f_i|_{[0,1]} \|_1\ .
\end{equation}
For $i\le n_0$,
$$|a_i| = |g_i(f)| \le |g_i (f-h)| + |g_i(h)|
< K\ep + \frac1{4n_0 \|f\|_1}\ .$$

For $i>n_0$, $|a_i| \le K(1+\ep)$.
Hence by \eqref{eq:fbbio}
\begin{align*}
1-\ep
&\le n_0 \left( K\ep + \frac1{4n_0\|f\|_1}\right) \|f\|_1
+ \sum_{i=n_0+1}^n K(1+\ep) \|f_i|_{[0,1]}\|_1 \\
&\qquad
< n_0 K\ep \|f\|_1 + \frac14 + \frac14 (1+\ep) < \frac34 < 1-\ep\ ,
\end{align*}
a contradiction, if $\ep< 1/4$.
\end{proof}

For $p=1$ we have a stronger result as a consequence of our next
theorem (Corollary~\ref{cor2.5}).

\begin{defn}
Let $1\le p<\infty$, $1/p + 1/q =1$.

a) Let $(f_i,g_i) \subseteq L_p(\real)
\times L_q (\real)$ be a frame for a subspace $X$ of $L_p(\real)$.
We say $(f_i,g_i)$ satisfies $(*)$ if
\begin{itemize}
\item[$(*)$]
for all $\ep>0$ and all bounded intervals $I\subseteq \real$, there
exists $n\in\nat$ so that for all $m\ge n$ and $f\in X$,
$$\Big\|\sum_{i=m+1}^\infty \langle f,g_i\rangle f_i|_I\Big\|_p
\le \ep \|f\|_p\ .$$
\end{itemize}

b)  A semi-normalized bounded minimal system $(f_i)_{i=1}^\infty$ in
$L_p(\real)$ satisfies $(**)$ if
\begin{itemize}
\item[$(**)$]
for all $\ep>0$ and bounded intervals $I\subseteq \real$ there exists
$n\in\nat$ so that for all $n<m\le m_1\le m_2$ and
$f= \sum_{i=1}^{m_2} a_i f_i\ \text{ with }\ \|f\|_p =1\ ,
\ \|\sum_{i=m}^{m_1} a_i f_i|_I \| \le \ep\ .$
\end{itemize}
\end{defn}

\begin{thm}\label{thm2.4}
Let $(f_i,g_i)_{i=1}^\infty$ be a frame or a semi-normalized bounded
fundamental minimal system for a subspace $X$ of $L_p(\real)$,
$1\le p<\infty$, satisfying $(*)$ or $(**)$, respectively.
Then $X$ embeds almost isometrically into $\ell_p$.
\end{thm}

$X$ embeds {\em almost isometrically into $\ell_p$} means that for
all $\ep >0$ there exits $T :X\to \ell_p$ with $(1+\ep)^{-1} \le \|Tf\|
\le 1+\ep$ for all $f\in S_X$.
The proof of Theorem~\ref{thm2.4} will yield, for all $\ep >0$, a partition
$\Pi = (D_s)_{s=1}^\infty$ of $\real$ into intervals so that for all
$f\in S_X$,
$$\|f - \EE_\Pi f\|_p <\ep\ .$$
$\EE_{\Pi}$ is the conditional expectation projection
$$f\longmapsto \sum_{s=1}^\infty \bigg( \int_{D_s} f\bigg)
\frac{\chix_{D_s}}{m(D_s)}\ .$$
Of course, in $L_p$, $(\frac{\chix_{D_s}}{m(D_s)})$ is 1-equivalent
to the unit vector basis of $\ell_p$.


{From} Proposition~\ref{prop2.1} and Theorem~\ref{thm2.4} we obtain

\begin{cor}\label{cor2.5}
If $(f_i,g_i)$ is a frame or a fundamental bounded minimal system
for a subspace $X$ of $L_1(\real)$ where the
$f_i$'s are uniformly discrete translates of some $f\in L_1(\real)$,
then $X$ embeds almost isometrically into $\ell_1$.
\end{cor}

\begin{proof}[Proof of Theorem~\ref{thm2.4}]
We first consider the frame case and let $C$ be the frame constant.
Thus for all $f\in X$ and $n\in \nat$,
$$\Big\| \sum_{i=1}^n \langle f,g_i\rangle f_i \Big\|_p \le C\|f\|_p\ .$$

Let $\ep>0$.
We inductively choose increasing sequences $(m_k)$ and $(n_k)$ in $\nat$
to obtain, where $I_k = [-m_k,m_k]$,
\begin{align}
&\text{for }\ f\in X,\ \text{ and }\ n\ge n_k\ ,\quad
\Big\| \sum_{i=n+1}^\infty \langle g_i,f\rangle f_i|_{I_{k-1}} \Big\|_p
\le \ep 2^{-k} \|f\|_p\ ,\hskip.5truein
\label{eq2.1}\\
&\text{for }\ f\in \text{span}\{f_i : i\le n_k\}\ ,\quad
\|f|_{\real\setminus I_k} \|_p
\le \ep 2^{-k} \|f\|_p\ .
\label{eq2.2}
\end{align}

We do this by setting $I_0 = \emptyset$, letting $n_1$ be arbitrary and
choose $m_1$ to satisfy \eqref{eq2.2} for $k=1$.
Then choose $n_2$ to satisfy \eqref{eq2.1} using $(*)$ and continue in
this manner.
We let $A_k = I_k \setminus I_{k-1}$, for $k\in\nat$.

Choose a partition $\pi_k$ of $A_k$ into intervals, $k\ge 1$, so that for all
$f\in \text{span}\{f_i :i\le n_{k+1}\}$,
\begin{equation}\label{eq2.3}
\Big\| f|_{A_k} - \sum_{D\in \pi_k} \frac{\chi_D}{m(D)}
\int_D f(x)\,dx \Big\|_p
\le \ep 2^{-k} \|f\|_p\ .
\end{equation}

Let $f\in X$ with $\|f\|_p =1$.
Then, with $n_0 =0$,
\begin{equation*}
\begin{split}
1& = \Big\|\sum_{i=1}^\infty \langle g_i,f\rangle f_i\Big\|_p
= \Big\| \sum_{s=1}^\infty \bigg( \sum_{i=n_{s-1}+1}^{n_s}
\langle g_i,f\rangle f_i|_{I_s}\bigg)
+ \sum_{s=1}^\infty \bigg( \sum_{i=n_{s-1}+1}^{n_s}
\langle g_i ,f\rangle f_i|_{\real\setminus I_s} \bigg) \Big\|_p\\
\noalign{\vskip6pt}
& \le \Big\| \sum_{s=1}^\infty \sum_{i=n_{s-1}+1}^{n_s}
\langle g_i,f\rangle f_i|_{I_s}\Big\|_p
+ 2C\ep\ ,\ \text{ by \eqref{eq2.2}}\\
\noalign{\vskip6pt}
&\le \Big\| \sum_{s=1}^\infty \sum_{i=n_{s-1}+1}^{n_s}
\langle g_i,f\rangle f_i|_{I_s\setminus I_{s-2}} \Big\|_p
+ 2C\ep + 2\ep\ ,\ \text{ by \eqref{eq2.1}}
\end{split}
\end{equation*}
where we let $I_{-1} = I_0 = \emptyset$
\begin{equation*}
= \Big\| \sum_{s=1}^\infty\ \sum_{i=n_{s-1}+1}^{n_s}
\langle g_i,f\rangle f_i |_{A_s\cup A_{s-1}} \Big\|_p
+ 2C\ep+2\ep \hskip.5truein
\end{equation*}
where we let $A_0 =\emptyset$
\begin{align}
&= \Big\| \sum_{s=1}^\infty \chi_{A_s}\ \sum_{i=n_{s-1}+1}^{n_{s+1}}
\langle g_i,f\rangle f_i\Big\|_p + 2C\ep +2\ep \notag\\
\noalign{\vskip6pt}
&\le \bigg( \sum_{s=1}^\infty \sum_{D\in \pi_s} \Big|
\int_D \sum_{i=n_{s-1}+1}^{n_{s+1}}
\langle g_i,f\rangle f_i(x)\, dx\Big|^p\bigg)^{1/p} + 4C\ep + 2\ep
\ \text{ by \eqref{eq2.3} .}
\label{eq2.4}
\end{align}

Now by \eqref{eq2.1} for $s\in\nat$,
\begin{equation}\label{eq2.5}
\bigg( \sum_{D\in\pi_s} \Big| \int_D \sum_{i=n_{s+1}+1}^\infty
\langle g_i,f\rangle f_i (x)\,dx \Big|^p\bigg)^{1/p}
\le \Big\| \sum_{i=n_{s+1}+1}^\infty \langle g_i,f\rangle f_i|_{I_s}\Big\|_p
\le \ep 2^{-(s+1)}\ .
\end{equation}

If $s>1$ then by \eqref{eq2.2},
\begin{equation}\label{eq2.6}
\bigg( \sum_{D\in \pi_s} \Big| \int_D \sum_{i=1}^{n_{s-1}}
\langle g_i,f\rangle f_i (x)\,dx\Big|^p\bigg)^{1/p}
\le 2^{-s+1} C\ep\ .
\end{equation}

{From} \eqref{eq2.4}, \eqref{eq2.5} and \eqref{eq2.6} we obtain that
\begin{equation*}
\begin{split}
1 = \|f\|_p
& \le \bigg( \sum_{D\in \bigcup\limits_{s=1}^\infty \pi_s}
\Big| \int_D f(x)\,dx \Big|^p\bigg)^{1/p}
+ \sum_{s=1}^\infty \ep 2^{-(s+1)}
+ \sum_{s=2}^\infty 2^{-s+1} C\ep + 6C\ep\\
\noalign{\vskip6pt}
& \le \bigg( \sum_{D\in \bigcup\limits_{s=1}^\infty \pi_s}
\Big| \int_D f(x)\,dx \Big|^p\bigg)^{1/p} + 8C\ep
= 1+8C\ep \ .
\end{split}
\end{equation*}
Thus $T : X\to \ell_p (\bigcup_{s=1}^\infty \pi_s)$ given by
$f\mapsto (\int_D f(x)\,dx)_{D\in \bigcup_{s=1}^\infty \pi_s}$
is the desired embedding.

The proof in the case that $(f_i)_{i=1}^\infty$ is a bounded fundamental
minimal system for $X\subseteq L_p(\real)$ is nearly identical.
We let $K = \sup_i \|g_i\|_q$ and in the construction replace
\eqref{eq2.1}--\eqref{eq2.3} by
\begin{align}
&\text{For all $n_k < n\le m\le \bar m$ and } f= \sum_1^{\bar m} a_i f_i
\in S_X\ ,
\tag{\ref{eq2.1}$'$}\\
&\qquad
\Big\| \sum_{i=n}^m a_i f_i|_{I_{k-1}}\Big\|_p \le \ep 2^{-k}\
\text{ (using $(**)$).}\notag\\
&\text{For all } \ f= \sum_{i=1}^{n_k} a_i f_i\ \text{ with $|a_i|\le K$
for $i\le n_k$, }\
\|f|_{\real\setminus I_k}\|_p \le \ep 2^{-k}\ .
\tag{\ref{eq2.2}$'$}\\
&\text{For all }\ f= \sum_{i=1}^{n_{k+1}} a_i f_i \ \text{ with $|a_i|\le K$
for $i\le n_{k+1}$, }
\tag{\ref{eq2.3}$'$}\\
&\qquad
\Big\| f|_{A_k} - \sum_{D\in\Pi_k} \frac{\chix_D}{m(D)} \int_D f(x)\,dx
\Big\|_p \le \ep 2^{-k}\ .\notag
\end{align}
The proof then proceeds as in the frame case for $f\in \text{span}(f_i)$,
$f= \sum_{i=1}^\ell a_i f_i$, $\|f\|_p =1$.
\end{proof}

\begin{remark}\label{shorterprf}
Let $X\subseteq L_p(\real)$ be as in Theorem~\ref{thm2.4} with $1<p<\infty$.
Then there is a shorter proof that yields $X\hookrightarrow \ell_p$.
In fact in the bounded minimal system case, one can replace $(**)$ by
the weaker
\begin{itemize}
\item[$(***)$]
For all $\ep>0$ and bounded intervals $I\subseteq \real$ there exists
$n\in \nat$ so that if\newline
$f\in \text{span}(f_i)_{i\ge n}$ with $\|f\|_p=1$,
then $\|f|_I\|_p < \ep$.
\end{itemize}

Indeed by \cite{KP}, \cite{J} and \cite{JO}, it suffices to prove
that if $(x_n)$ is a normalized weakly null sequence in $X$ then some
subsequence is 2-equivalent to the unit vector basis of $\ell_p$.
Then, from $(*)$ or $(***)$, it is easy to find $(x_{n_i})$ and intervals
$I_1 \subseteq I_2 \subseteq \cdots$ so that
$\|x_{n_i}|_{I_i\setminus I_{i-1}}\|_p > 1-\frac{\ep}{2^i}$ for all $i$
and deduce the result.
\end{remark}

We will say that a {\em frame $(f_i,g_i)_{i=1}^\infty$ for $X$
satisfies a lower $\ell_q$-estimate\/}
if for some $K<\infty$ and all $x\in X$,
\begin{equation*}
\bigg( \sum_{i=1}^\infty |g_i(x)|^q\bigg)^{1/q}
\le K\Big\| \sum_{i=1}^\infty g_i (x) f_i\Big\|
= K\|x\|\ .
\end{equation*}
A Hilbert frame, by definition, satisfies lower (and upper)
$\ell_2$-estimates.

If $(x_i)_{i=1}^\infty$ is a fundamental bounded minimal system for $X$
we say that $(x_i)_{i=1}^\infty$ satisfies a lower $\ell_q$-estimate
if for some $K$ and all scalars $(a_i)_1^n$
\begin{equation*}
\bigg(\sum_{i=1}^n |a_i|^q\bigg)^{1/q}
\le K \Big\| \sum_{i=1}^n a_i x_i\Big\|\ .
\end{equation*}

\begin{prop}\label{prop2.6}
Let $1<p<\infty$, $1/p + 1/q =1$.

{\em a)} Assume $1<p\le 2$ and let $(f_i)_{i=1}^\infty$ be a sequence of
uniformly discrete translations of $f\in L_p(\real)$.
Let either $(f_i,g_i)_{i=1}^\infty \subseteq L_p(\real)\times L_q(\real)$
be a frame or $(f_i)_{i=1}^\infty$ be a fundamental bounded minimal system
for $X\subseteq L_p(\real)$.
If $(f_i)_{i=1}^\infty$ admits a lower $\ell_q$-estimate, then $X$ embeds
almost isometrically into $\ell_p$.

{\em b)} Let $(f_i,g_i)_{i=1}^\infty$ be a frame for $L_p(\real)$, where
$(f_i)_{i=1}^\infty$ is a sequence of uniformly discrete translations of
$f\in L_p(\real)$.
Then for all bounded measurable sets $B$ of positive measure,
$\sum_{i=1}^\infty \|g_i|_B\|_1^q = \infty$.
\end{prop}

\begin{remark}\label{rem:hypothesis}
The hypothesis in a) would be vacuous for $p>2$ since some subsequence of
$(f_i)$ is equivalent to the unit vector basis of $\ell_p$.
\end{remark}


\begin{proof}
First let $(f_i,g_i)_{i=1}^\infty$ be a frame for $X\subseteq L_p(\real)$
as in a).
Assume for all $f\in X$,
$$\bigg( \sum_{i=1}^\infty |g_i(f)|^q\bigg)^{1/q} \le K\|f\|_p\ .$$
For any bounded interval $I\subseteq \real$, $f\in X$ and $n\in \nat$,
\begin{align*}
\Big\|\sum_{i=n}^\infty g_i (f) f_i\big|_I \Big\|_p
&
\le \sum_{i=n}^\infty |g_i (f)|\ \|f_i\big|_I\|_p\\
\noalign{\vskip6pt}
&
\le \bigg( \sum_{i=n}^\infty |g_i(f)|^q\bigg)^{1/q}
\bigg( \sum_{i=n}^\infty \|f_i\big|_I\|_p^p \bigg)^{1/p}\\
\noalign{\vskip6pt}
&
\le K \|f\|_p \bigg( \sum_{i=n}^\infty \|f_i\big|_I \|_p^p\bigg)^{1/p}\ .
\end{align*}
{From} Proposition~\ref{prop2.1} we obtain that $(*)$ holds and so
Theorem~\ref{thm2.4} applies.

Similarly, if $(f_i)$ is a fundamental bounded minimal system for $X$ and
$f= \sum_1^{\bar m} a_i f_i$ with $\|f\|_p=1$, we have
$$\Big\| \sum_n^m a_i f_i|_I\Big\|_p \le K\bigg( \sum_n^m \|f_i|_I\|_p
\bigg)^{1/p}$$
and so, again, we have $(**)$ and apply Theorem~\ref{thm2.4}.

b) Assume that for some $B$ of positive measure $\sum_{i=1}^\infty
\|g_i|_B\|_1^q <\infty$.
Let $h\in L_\infty (B)$, $|h| =1$.
So $h = \sum_{i=1}^\infty \langle h,g_i\rangle f_i|_B$, the series
converging in $L_1(B)$.
Thus
\begin{align*}
m(B) = \|h\|_1 & \le \sum_{i=1}^n |\langle h,g_i\rangle| \|f_i\|_1
+ \sum_{i=n+1}^\infty |\langle h,g_i|_B\rangle|\, \|f_i|_B\|_1\\
&\le \sum_{i=1}^n |\langle h,g_i\rangle|\, \|f_i\|_1
+ \bigg( \sum_{i=n+1}^\infty \|g_i|_B\|_1^q\bigg)^{1/q}
\bigg( \sum_{i=n+1}^\infty \|f_i|_B\|_1^p\bigg)^{1/p}\ .
\end{align*}
Then, as in the proof of Proposition~\ref{prop2.3bis}, we can choose $n$
so that the second term does not exceed $m(B)/4$, and given this $n$,
choose $h$ to make the first term also less than $m(B)/4$.
Thus $m(B) < \frac12 m(B)$, a contradiction.
\end{proof}

Part a) of Proposition~\ref{prop2.6} yields a quantitative improvement
of Theorem~\ref{thm:OZ}.
If ${\{\flambda :\lambda\in\Lambda\}}$ is unconditional basic in $L_2(\real)$,
then given $\ep>0$ there is a partition $\Pi$ of $\real$ so that for all
$g\in X_2 (f,\Lambda)$,
$$\|g- \EE_{\Pi} g\|_2 \le \ep \|g\|_2\ .$$

\begin{remark}\label{rem2.7}
a) Let $f= \chix_{[0,1]} - \chix_{[1,2]}\in L_p(\real)$. The
sequence $(f_{(n)})_{n\in\zed}$ is basic in $L_1(\real)$
when ordered as $(f_{(0)},f_{(1)},f_{(-1)},f_{(2)},
f_{(-2)},\ldots)$.
It is not unconditional since
\begin{equation*}
\Big\| \sum_{i=-2n}^{2n} f_{(i)}\Big\|_p = 2^{1/p}
\ \text{ but }\
\Big\| \sum_{i=-n}^n f_{(2i)}\Big\|_p = (2n+2)^{1/p}\ .
\end{equation*}
For ${1<p<\infty}$, $(f_{(n)})_{n\in\zed}$ is not a frame nor a minimal
system in $L_p(\real)$.
The latter follows easily from the fact that
$$\lim_{n\to\infty} \Big\| f_{(0)} + \sum_{k=1}^n \frac{n-k}n
(f_{(k)} + f_{(-k)})\Big\|_p =0\ .$$

b) (due to S.J.~Dilworth)
Let $1\le p <\infty$ and let
$$f = \chix_{[-3/2,\, -1/2]} + 2\chix_{[-1/2,\, 1/2]}
+\chix_{[1/2,\, 3/2]}\ .$$
For $n\in \nat$ set
$$g_n = f + \sum_{k=1}^n (-1)^n \ \frac{n-k+1}{n}\
(f_{(k)} + f_{(-k)})\ .$$
Then for $x\ge 0$
\begin{equation*}
g_n(x) =
\begin{cases}
\ds -f_{(-1)} (x) + f_{(0)}(x) - f_{(1)}(x) = -1 +2-1 = 0
&\text{if  $x\in [0,\frac12]$}\\
\noalign{\vskip6pt}
 f_{(0)}(x\!-\! f_{(1)}(x)\!+\!\frac{n-1}n f_{(2)} (x)
= 1\!-\!2\!+\!\frac{n-1}n\! =\! -\! \frac1n
&\text{if $x \in [\frac12, \frac32]$}\\
\noalign{\vskip6pt}
\ds (-1)^k \frac{n-k+1}n f_{(k)}(x) + (-1)^{k+1} \frac{n-k}n f_{(k+1)}(x) &\\
\hskip1in
+ (-1)^{k+2} \frac{n-k-1}n f_{(k+2)} (x) =0
 &\text{if  $x\in [ k\!+\!\frac12, k\!+\!\frac32]$}\\
\noalign{\vskip4pt}
&\text{for some $1\!\le\!k\!\le\! n\!-\!2$}\\
\noalign{\vskip6pt}
\ds (-1)^{n-1}\frac2n f_{(n-1)}(x) + (-1)^n \frac1n f_{(n)}(x)=0
&\text{if  $x\in [n \!-\!\frac12, n\!+\!\frac12]$}\\
\noalign{\vskip6pt}
\ds (-1)^n \frac1n f_{(n)} (x) = (-1)^n \frac1n
&\text{if  $x\in [n\!+\!\frac12, n\!+\!\frac32]$}\ .
\end{cases}
\end{equation*}
Thus $\|g_n\|_p = 4^{1/p}/n$, hence $f_{(0)}\in [\{f_{(k)} :k\in\zed\setminus
\{0\}\}]$ and so $(f_{(k)})_{k\in\zed}$ is not a minimal system in
$L_p(\real)$.
Furthermore,
$$\hat\chix_{[-1/2,\, 1/2]} (x)
= \frac1{\sqrt{2\pi}} \int_{-1/2}^{1/2} e^{-ixt} \,dt
= \frac1{\sqrt{2\pi}}\ \frac{\sin x}{x}\ .$$
It follows that
$$\hat f (x) = \frac1{\sqrt{2\pi}}\ \frac{\sin x}{x} [2+\cos x]$$
so $\hat f(x) \ne0$ a.e..

c) (due to D.~Freeman)
It is well known that if $(f_i)_{i=1}^\infty$ is a normalized sequence
in $L_1(\real)$ with $\|f_i|_{I_i}\| \ge \lambda > \frac12$ for all $i$
and some sequence of pairwise disjoint measurable sets $I_i\subseteq \real$,
then $(f_i)_{i=1}^\infty$ is equivalent to the unit vector basis of $\ell_1$.
Indeed
\begin{align*}
\| \sum a_i f_i\|_1
& \ge \| \sum a_i f_i|_{I_i}\|_1
- \sum |a_i|\ \|f_i|_{\real_i\setminus I_i} \|_1\\
& \ge \lambda \sum |a_i| - (1-\lambda) \sum |a_i|
= (2\lambda -1) \sum |a_i|\ .
\end{align*}
Thus if $\{\flambda :\lambda \in\Lambda\}$ is a sequence of uniformly
discrete translations of $0\ne f\in L_1(\real)$, then it can be split into
a finite number of subsequences, each equivalent to the unit vector basis
of $\ell_1$.

d) 
By Theorem~\ref{thm:ER},
if $1\le p< \infty$ and $f\in L_p(\real)$, $f\ne 0$ then $\{\flambda :\lambda
\in\real\}$ is linearly independent (see also Theorem~\ref{thm4.6} below).
However one can find $f\in L_1(\real)$ so that $\{f_{(n)} :n\in\zed\}$ is
not $\omega$-linearly independent in its natural order \cite{R}.
\end{remark}

We next turn to the case where $(f_i)$ is unconditional basic in $L_p$.
We first recall

\begin{prop}\label{prop2.9}
\cite[Lemma 2]{JO}.
Let $1\le p\le 2$.
Let $(f_i) \subseteq L_p (\real)$ be seminormalized and unconditional basic.
Assume that for some $\delta >0$ there exists a sequence of disjoint
measurable sets $(B_i)_{i=1}^\infty$ with $\|f_i|_{B_i}\|_p \ge \delta$,
for all $i$.
Then $(f_i)_{i=1}^\infty$ is equivalent to the unit vector basis of $\ell_p$.
\end{prop}

\begin{cor}\label{cor2.10}
Let $(f_i)_{i=1}^\infty$ be an unconditional basic sequence in
$L_p(\real)$, $1\le p\le2$.
Assume the $f_i$'s are all translates of some fixed $f\in L_p(\real)$.
Then $(f_i)_{i=1}^\infty$ is equivalent to the unit vector basis of $\ell_p$.
\end{cor}

\begin{proof}
Let $f_i = f_{(\lambda_i)}$ for $i\in\nat$.
Let $\rho \equiv \frac12 \inf \{|\lambda_i - \lambda_j| :i\ne j\} >0$.
Let $I$ be an interval of length $\rho$ with $\|f|_I\|_p = \delta >0$.
If $B_i = I+\lambda_i$ for $i\in \nat$ then the $B_i$'s are pairwise
disjoint and $\|f_i|_{B_i} \| = \| f|_I\| = \delta$ for all $i$.
Proposition~\ref{prop2.9} yields the result.
\end{proof}

As we shall see the situation is more complicated for $p>2$, and
especially so for $p>4$.

\begin{thm}\label{thm2.11}
Let $2<p\le 4$ and let $(f_i)\subseteq L_p (\real)$ be an unconditional
basis for $X \subseteq L_p(\real)$.
Assume the $f_i$'s are all translates of some fixed $f\in L_p(\real)$.
Then $X$ embeds isomorphically into $\ell_p$.
\end{thm}

\begin{lem}\label{lem2.12}
Let $p\ne 2$ and let $X$ be a subspace of $L_p(\real)$ not containing
an isomorph of $\ell_p$.
Then there exists $c>0$ so that $\NORM{f} = \|f|_{[-c,c]}\|_p$
is an equivalent norm on $X$.
\end{lem}

\begin{proof}
If the lemma is false then we can find $(f_k)_{k=1}^\infty \subseteq S_X$
and $(m_k)_{k=1}^\infty \subseteq \nat$ so that
$\| f_k|_{[-m_k,m_k]} \| \ge 1-2^{-2k-1}$ and
$\| f_{k+1} |_{[-m_k,m_k]} \| \le 2^{-2k-1}$ for all $k\in \nat$.
It follows easily that $(f_k)_{k=1}^\infty$ is equivalent to
$(f_k |_{[-m_k,m_k]\setminus [-m_{k-1},m_{k-1}]} )_{k=1}^\infty$
which, being seminormalized and disjointly supported, is equivalent
to the unit vector basis of $\ell_p$.
\end{proof}

We shall also use

\begin{prop}\label{prop2.13} \cite{JO}.
Let $X$ be a subspace of $L_p(\real)$, $2<p<\infty$, which does not
contain an isomorph of $\ell_2$.
Then $X$ embeds isomorphically into $\ell_p$.
\end{prop}
\noindent
In fact by \cite{KW}, $X$ must then embed almost isometrically into $\ell_p$.

We set some notation and recall some things before proving the theorem.
We let $(h_i)$ denote the normalized Haar basis for $L_p[0,1]$ regarded,
canonically, as a subspace of $L_p(\real)$.
As mentioned in the introduction,
for $i\in \nat$ and $n\in\zed$ we let $h_{(i,n)} (\cdot) = h_i((\cdot)-n)$.
Thus, $(h_{(i,n)})_{i\in \nat,\ n\in \zed}$
is an unconditional basis for $L_p(\real)$.

G.~Schechtman \cite{S} made the very useful observation that if
$(f_i)_{i=1}^\infty$ and $(g_i)_{i=1}^\infty$ are seminormalized
unconditional basic sequences in $L_p(\real)$, $1<p<\infty$, with
\begin{equation}\label{eq2.8}
\sum_{i=1}^\infty \big\| \, |f_i| - |g_i |\,\big\|_p < \infty
\end{equation}
then$(f_i)_{i=1}^\infty$ is equivalent to $(g_i)_{i=1}^\infty$.
This follows from (1.3).    
In particular, if $(f_i)_{i=1}^\infty$ is seminormalized unconditional
basic in $L_p(\real)$ then, by first approximating each $(f_i)_{i=1}^\infty$
by a simple dyadic function and then using the above consequence of
(1.3), 
there exists a block basis $(g_i)_{i=1}^\infty$  of
$(h_{(i,n)})_{i\in\nat,\, n\in\zed}$ satisfying \eqref{eq2.8} and thus
being equivalent to $(f_i)_{i=1}^\infty$.

\begin{proof}[Proof of Theorem~\ref{thm2.11}]
By Proposition~\ref{prop2.13}, it suffices to prove that $X$ does not
contain an isomorph of $\ell_2$.
By our above remarks we can choose a block basis $(g_i)_{i=1}^\infty$
of $(h_{(i,n)})$ which satisfies \eqref{eq2.8}.
In particular $(g_i)_{i=1}^\infty$ is equivalent to $(f_i)_{i=1}^\infty$
and we maintain
\begin{equation}\label{eq2.9}
\sum_{i=1}^\infty \|g_i|_I \|_p^p < \infty\ \text{ for all bounded
intervals }\ I.
\end{equation}

Thus we need only show that $[(g_i)_{i=1}^\infty]$ does not contain an isomorph
of $\ell_2$.
If this is false, then there exists a normalized block basis
$(\bar g_i)_{i=1}^\infty$ of $(g_i)_{i=1}^\infty$
which is equivalent to the unit vector basis of $\ell_2$.
Set $\bar X = [(\bar g_i)_{i=1}^\infty]$.
By Lemma~\ref{lem2.12} there exists $M\in\nat$ and
$1\le C< \infty$ so that for all $\bar g \in \bar X$, $I = [-M,M]$,
\begin{equation}\label{eq2.10}
\|\bar g|_I \|_p \ge C^{-1} \|\bar g\|_p\ .
\end{equation}

Since $(g_i)_{i=1}^\infty$ is a block basis of $(h_{(i,n)})$ then so is
the normalized and unconditional sequence
$(g_i|_I /\|g_i |_I \|_p)_{i=1}^\infty$.
This yields lower $\ell_p$ and upper $\ell_2$ estimates for this sequence
and $(g_i)$ (see (1.2)).
{From} this and \eqref{eq2.10} we obtain for some constant $D<\infty$
and for all $(a_i)\subseteq \real$, and
$\bar g = \sum_{i=1}^\infty a_i g_i \in \bar X$,
\begin{equation}\label{eq2.11}
\begin{split}
D^{-1}\bigg( \sum_{i=1}^\infty |a_i|^p\bigg)^{1/p}
&\le D^{-1/2} \Big\| \sum_{i=1}^\infty a_i g_i\Big\|_p\\
\noalign{\vskip6pt}
&\le \Big\| \sum_{i=1}^\infty a_i g_i|_I \Big\|_p
= \Big\| \sum_{i=1}^\infty a_i \| g_i|_I \|_p
\frac{g_i|_I}{\|g_i|_I\|_p} \Big\|_p\\
\noalign{\vskip6pt}
&\le D\bigg( \sum_{i=1}^\infty |a_i|^2 \|g_i|_I \|_p^2\bigg)^{1/2}\ .
\end{split}
\end{equation}

By \eqref{eq2.9} there exists $n_0 \in \nat$ with
\begin{equation}\label{eq2.12}
\bigg( \sum_{i=n_0}^\infty \| g_i |_I\|_p^p\bigg)^{1/p}  < (2C)^{-1} D^{-2}\ .
\end{equation}

Let $\bar g$ be an element of $S_X$ which has the property that if we
expand it in terms of the $g_i$'s, i.e., if we write it as
$\bar g = \sum_{i=1}^\infty a_i g_i$ then $a_j=0$ for $j\le n_0$.
{From} \eqref{eq2.10} and \eqref{eq2.11},
\begin{equation*}
C^{-1} \le \|\bar g|_I\|_p
\le D\bigg( \sum_{i=n_0}^\infty a_i^2 \|g_i|_I \|_p^2\bigg)^{1/2}
\le D\Big[ \| (a_i^2)_{i=n_0}^\infty \|_{\ell_{p/2}}
\cdot \|\big( \|g_i|_I\|_p^2\big)_{i=n_0}^\infty \|_{\ell_{\frac{p}{p-2}}}
\Big]^{1/2}
\end{equation*}
(applying H\"older's inequality for $p/2$ and $p/p-2$)
\begin{equation*}
\begin{split}
&= D \|(a_i)_{i=n_0}^\infty \|_{\ell_p}
\bigg( \sum_{i=n_0}^\infty \|g_i |_I \|_p^{\frac{2p}{p-2}}
\bigg)^{\frac{p-2}{2p}}\\
\noalign{\vskip6pt}
& \le D^{2} \| \bar g|_I\|_p
\bigg(\sum_{i=n_0}^\infty \|g_i |_I \|_p^p\bigg)^{1/p} \\
\end{split}
\end{equation*}
(by \eqref{eq2.11} and since $p\le 4$, $\frac{2p}{p-2}\ge p$)
$$\le (2C)^{-1}\quad \text{by \eqref{eq2.12},}\hskip.9truein $$
which is a contradiction.
\end{proof}

When $p>4$ the possible structure is more complicated.

\begin{thm}\label{thm2.13}
Let $4< p<\infty $.
There exists $f\in L_p(\real)$ and $\Lambda \subseteq \zed$ so that
$(f_{(\lambda)})_{\lambda\in\Lambda}$ is an unconditional basic sequence with
$X_p (f,\Lambda)$ containing an isomorph of $L_p(\real)$.
\end{thm}

\begin{proof}
We identify, in the canonical way, $L_p(\real)$ with $(\bigoplus_{i\in\zed}
L_p [0,1])_{\ell_p}$.
Since $L_p[0,1]$ is isometrically isomorphic to $L_p([0,1]^2)$, we need
only produce
$f= (f_i)_{i\in \zed} \in (\bigoplus_{i\in\zed} L_p[0,1]^2)_{\ell_p}$
and $\Lambda\subseteq \nat$ so that setting for $\lambda\in\Lambda$
$f_{(\lambda)} = (f_{i-\lambda})_{i\in\zed}$, then
$X_p (f, \Lambda)$ contains an isomorph of $L_p[0,1]$
and $(f_{(\lambda)})_{\lambda\in\Lambda}$ is unconditional.

Letting, as before, $(h_n)_{n=1}^\infty$ be the normalized Haar basis
for $L_p[0,1]$ and $(r_n)_{n=1}^\infty$ the Rademacher functions on
$[0,1]$ we have, for some constants $C_p$ and $D_p$ (see (1.3)),
for all $(a_i)\subseteq \real$,
\begin{equation}\label{eq2.13}
\bigg(\sum_{i=1}^\infty |a_i|^2\bigg)^{1/2}
\le \Big\| \sum_{i=1}^\infty a_i r_i\Big\|_p
\le C_p \bigg( \sum_{i=1}^\infty |a_i|^2\bigg)^{1/2}
\end{equation}
and
\begin{equation}\label{eq2.14}
\Big\| \sum_{i=1}^\infty a_i h_i \Big\|_p
\ \buildrel {D_p}\over\sim \
\Big\| \sum_{i=1}^\infty a_i^2 |h_i|^2\Big\|_{p/2}^{1/2}\ .
\end{equation}

Since $p>4$ we can choose $(\ep_i)_{i=1}^\infty \subseteq (0,1)$ so that
\begin{equation}\label{eq2.15}
\sum_{i=1}^\infty \ep_i^p =1
\end{equation}
and there exists a partition $(J_n)_{n=1}^\infty$ of $\nat$ into finite
intervals with
\begin{equation}\label{eq2.16}
\sum_{j\in J_n} \ep_j^4 =1\ \text{ for all }\ n\in \nat\ .
\end{equation}

We are ready to define $f= (f_i)_{i\in\zed} \in (\bigoplus_{i\in\zed}
L_p [0,1]^2)_{\ell_p}$.
Set for $i\in \zed$,
\begin{equation}\label{eq2.17}
f_i = \begin{cases}
\ep_j h_n \otimes r_j\ ,\ \text{ if }\ i=3^j\ \text{ with }\ j\in J_n
\ \text{ for some }\ n\in\nat\\
\noalign{\vskip6pt}
0\ ,\ \text{ otherwise.}
\end{cases}
\end{equation}
where $h_n\otimes r_j$ is placed on the $i^{th}$ copy of $[0,1]^2$.
Note that
$$\|f\|_p^p = \sum_{n\in\nat} \sum_{j\in J_n} \|\ep_j h_n \otimes r_j\|_p^p
= \sum_{n\in\nat} \sum_{j\in J_n} \ep_j^p =1\ .$$
Let $\Lambda = \{ -3^j :j\in \nat\}$ and so our translated sequence is
$(f_{(-3^j)})_{j=1}^\infty$.
For ease of notation below we shall write $f_{(-3^j)}$, $f$ shifted $3^j$
units left, as $f^{(-3^j)}$, and $f^{(-3^j)} = (f_i^{(-3^j)})_{i\in\zed}$
where $f_i^{(-3^j)}$ denotes $f^{(-3^j)}$ restricted to the $i^{th}$ $[0,1]^2$.

Now $f_0^{(-3^j)} = \ep_j h_n \otimes r_j$ if $j\in J_n$ and so for
$(a_j) \subseteq \real$,
\begin{equation}\label{eq2.18}
\begin{split}
\Big\| \sum_{j\in\nat} a_j f_0^{(-3^j)}\Big\|_p^p
& = \Big\| \sum_{n\in \nat} \sum_{j\in J_n} a_j \ep_j h_n\otimes r_j
\Big\|_p^p \\
\noalign{\vskip6pt}
& = \int_0^1 \int_0^1 \Big| \sum_{n\in \nat} \sum_{j\in J_n} a_j \ep_j
h_n(s) r_j (t)\Big|^p\, dt\,ds\\
\noalign{\vskip6pt}
& \buildrel {C_p^p}\over \sim\
\int_0^1 \Big|\sum_{n\in\nat}\sum_{j\in J_n} a_j^2 \ep_j^2 h_n^2 (s)\Big|^{p/2}
\, ds\ ,\ \text{by \eqref{eq2.13}}\\
\noalign{\vskip6pt}
&= \Big\| \sum_{n\in \nat} \bigg( \sum_{j\in J_n} a_j^2 \ep_j^2\bigg) h_n^2
\Big\|_{p/2}^{p/2} \\
\noalign{\vskip6pt}
& \buildrel {D_p^p}\over \sim\
\Big\| \sum_{n\in \nat} \bigg( \sum_{j\in J_n} a_j^2 \ep_j^2\bigg)^{1/2}
h_n \Big\|_p^p\ ,\ \text{ by \eqref{eq2.14}.}
\end{split}
\end{equation}
Now for $j\in \nat$, $f_\ell^{(-3^j)}\ne0$ iff $\ell = 3^k - 3^j$ for some
$k\in \nat$.
If $\ell\ne0$ and $\ell = 3^k - 3^j = 3^{k'} - 3^{j'}$ for
$k,k',j,j'\in\nat$ then $k= k'$ and $j= j'$.
Thus the functions $(f^{(-3^j)})_{j\in \nat}$ are
disjointly supported except on the $0^{th}$ copy of $[0,1]^2$.
Also
$$\Big\| \sum_{\substack{\ell\ne0\\ \ell\in\zed}} f_\ell^{(-3^j)}\Big\|_p^p
= 1-\ep_j^p\ ,$$
{From} this and \eqref{eq2.18} we obtain for some $K$, for all $(a_i)\subseteq
\real$,
\begin{equation}\label{eq2.19}
\Big\| \sum_{j\in \nat} a_j f^{(-3^j)}\Big\|_p^p
\ \buildrel K\over\sim\
\Big\| \sum_{n\in\nat} \bigg( \sum_{j\in J_n} a_j^2 \ep_j^2\bigg)^{1/2}
h_n \Big\|_p^p
+ \sum_{j\in \nat} |a_j|^p\ .
\end{equation}
Thus $(f^{(-3^j)})_{j=1}^\infty$ is unconditional and we shall next construct
a block basis $(b^{(n)})_{n=1}^\infty$ of $(f^{(-3^j)})_{j=1}^\infty$
which is equivalent to $(h_n)_{n=1}^\infty$.

For $n\in\nat$ set
$$b^{(n)} = \sum_{j\in J_n} \ep_j f^{(-3^j)}\ .$$
{From} \eqref{eq2.19},  for $(c_n) \subseteq \real$
\begin{equation*}
\begin{split}
\Big\| \sum_{n=1}^\infty c_n b^{(n)}\Big\|_p^p
&=  \Big\| \sum_{n=1}^\infty c_n \sum_{j\in J_n} \ep_j f^{(-3^j)}\Big\|_p^p\\
&\buildrel {K^p} \over\sim\ \Big\| \sum_{n=1}^\infty c_n
\bigg( \sum_{j\in J_n} \ep_j^4\bigg)^{1/2} h_n \Big\|_p^p
+ \sum_{n=1}^\infty \sum_{j\in J_n} |c_n|^p \ep_j^p\\
& = \Big\| \sum_{n=1}^\infty c_n h_n \Big\|_p^p
+ \sum_{n=1}^\infty |c_n|^p \bigg( \sum_{j\in J_n} \ep_j^p\bigg)\ .
\end{split}
\end{equation*}
Thus, using this and (1.2), the lower $\ell_p$ estimate of
$(h_n)_{n=1}^\infty$, we see that $(b^{(n)})_{n=1}^\infty$ is equivalent
to $(h_n)_{n=1}^\infty$.
\end{proof}

We next note that under certain additional assumptions we cannot have the
situation  of Theorem~\ref{thm2.13}.

\begin{prop}\label{prop2.14}
Let $4 < p<\infty$ and let $(f_i)_{i=1}^\infty$ be an unconditional basis
for $X \subseteq L_p (\real)$ where the $f_i$'s are all translations of
some fixed $f\in L_p(\real)$.
If either
\begin{itemize}
\item[a)] $f\in L_2 (\real)$ or
\item[b)] $\sum_{n\in \zed} \|f |_{[n-1,n]}
\|_{_{_{\scriptstyle p}}}^{\frac{2p}{p-2}} <\infty$
\end{itemize}
then $X$ embeds isomorphically into $\ell_p$.
\end{prop}

\begin{proof}
b) follows easily from the proof of Theorem~\ref{thm2.11}.
Indeed we can use b) to deduce the next to last inequality in that proof,
rather than using $p\le4$ as was done there.

a) We assume the contrary so by Proposition~\ref{prop2.13} $X$ contains an
isomorph of $\ell_2$.
We choose $I$, $(g_i)_{i=1}^\infty$ and $(\bar g_i)_{i=1}^\infty$
as in the proof of Theorem~\ref{thm2.11} with the additional assumption
that $\sum_{i=1}^\infty \big\| \, |f_i | - |g_i|\,\big\|_2 <\infty$.

Hence, using $f\in L_2(\real)$,
\begin{equation}\label{eq2.20}
\sum_{i=1}^\infty \| g_i |_I \|_2^2 < \infty\ .
\end{equation}

Now $(\bar g_i|_I)_{i=1}^\infty$ is a block basis of $(h_{(i,n)})$ which
is equivalent to the unit vector basis of $\ell_2$.
This forces $\|\cdot \|_p$ and $\|\cdot\|_2$ to be equivalent on
$[(\bar g_i|_I)_{i=1}^\infty]\subseteq L_p(I)$ \cite{KP}.
$(\bar g_i)_{i=1}^\infty$ is also a normalized block basis of
$(g_i)_{i=1}^\infty$ and so we may write
$\bar g_i = \sum_{j= n_{i-1}+1}^{n_i} c_j g_j$
for some scalars $(c_j)$, $n_0< n_1 < \cdots$ and all $i\in\nat$.
Since $(g_i|_I)_{i=1}^\infty$ is also a block basis of
$(h_{(i,n)})$ and hence is orthogonal in $L_2(I)$, we have for $i\in\nat$,
\begin{equation*}
\begin{split}
\|\bar g_i |_I\|_2
& = \Big\| \sum_{j= n_{i-1}+1}^{n_i} c_j g_i |_I\Big\|_2
= \bigg( \sum_{j=n_{i-1}+1}^{n_i} c_j^2 \|g_j |_I \|_2^2\bigg)^{1/2}
 \le \sup_j|c_j| \bigg( \sum_{j=n_{i-1}+1}^{n_i} \|g_i |_I \|_2^2\bigg)^{1/2}
\end{split}
\end{equation*}
and the latter converges to $0$ as $i\to\infty$ by \eqref{eq2.20}.
Thus $\|\bar g_i |_I \|_2 \to 0$ so $\|\bar g_i|_I\|_p \to 0$
which is a contradiction.
\end{proof}

We next present two more examples.
The first is an easy example of a translation sequence in $L_p$
($2<p$) which is unconditional but not equivalent to the $\ell_p$ basis and
so Theorem~\ref{thm2.11} cannot be improved to get $(f_i)$ equivalent to
the unit vector basis of $\ell_p$.
The second is a translation sequence $(f_i)$ in $L_p$, $p>4$, which is
basic but not unconditional.

\begin{example}\label{ex2.16}
Let $2<p<\infty$.
There exists $f\in L_p(\real)$ so that $(f_{(n)})_{n=1}^\infty$, the
sequence of translations of $f$ by $n\in\nat$, is unconditional basic
but not equivalent to the unit vector basis of $\ell_p$.

Of course we already know this for $p>4$ by Theorem~\ref{thm2.13}.
Let $(r_n)_{n\in \zed}$ be an enumeration of
the Rademacher functions on $[0,1]$ extended
trivially to functions defined on all of $\real$.
We define $\tilde r_n(\cdot) = r_n ((\cdot) -n)$, for $n\in\zed$, and
let $f= \sum_{n\in\zed} \frac{\tilde r_n}{\sqrt{|n|}}$,
where we regard $\frac1{\sqrt{|0|}} = 1$.
Note that $\|f\|_p^p = 1+2\sum_{n=1}^\infty n^{-p/2} <\infty$, since $p>2$.
For $(a_i) \in c_{00}$, $g= \sum a_i f_{(i)}$ and $x\in [k,k+1]$,
$k\in\zed$, we observe
$$g(x) = \sum_{i=1}^\infty a_i f (x-i)
= \sum_{i=1}^\infty a_i \frac{r_{k-i} (x-k)}{\sqrt{|k-i|}}\ .$$
Thus, for some $c_p>0$
$$\|g\|_p^p = \sum_{k\in\zed} \| g|_{[k,k+1]}\|_p^p
\ \buildrel {c_p}\over \sim\
\sum_{k\in\zed} \bigg( \sum_{i=1}^\infty \frac{a_i^2}{|k-i|}\bigg)^{p/2}\ ,$$
which shows that $(f_{(i)})_{i=1}^\infty$ is unconditional.
Moreover, if we let $a_i=1$, for $i=1,\ldots,m\in\nat$ for $m\in\nat$ we
obtain
$$\Big\| \sum_{i=1}^m f_{(i)}\Big\|_p^p
\ge c_p \sum_{k=1}^m
\bigg( \sum_{i=1}^m \frac1{|m-i|}\bigg)^{p/2}
\ge c_p\, m(\log m)^{p/2}\ .$$
Thus $(f_{(i)})$ is not equivalent to the unit vector basis
of $\ell_p$.\qed
\end{example}
\medskip

\begin{example}\label{ex2.17}
Let $p>4$.
There exists $f\in L_p(\real)$ and $\Lambda\subset \zed$ so that
$\{f_{(\lambda)} :\lambda\in\Lambda\}$ is basic in  some order, but
not unconditional.

As in Theorem~\ref{thm2.13} we identify $L_p(\real)$ with $(\oplus_{n\in\zed}
L_p[0,1])_p$, and we write $f$ as $(f_i:\in\zed)$ with $f_i\in L_p(0,1)$,
for $i\in \zed$, and, as in Theorem~\ref{thm2.13}, we write $f^{(\lambda)}$
instead of $f_{(\lambda)}$.

For $j\in \nat$ let $a_j = j^{-1/4}$, and $a_0=1$.
Let $(r_j)$ be the Rademacher sequence on $[0,1]$.
We define $f= (f_i)_{i\in\zed}$ by
$$f_i = \begin{cases}
a_{j-1} r_j - a_{j+1} r_{j+1}&\text{if $i=3^j$ for some $j\in \nat$,}\\
\noalign{\vskip6pt}
0&\text{otherwise.}
\end{cases}$$
Since $p>4$, $(a_i) \in \ell_p$ and, thus,
$f\in (\oplus_{n\in\zed} L_p [0,1])_p$.
We let $\Lambda = \{-3^j :j\in \nat\}$.
For $(b_j)_{j=1}^n \subset \real$ we compute $(b_0=0)$
\begin{equation*}
\begin{split}
\sum_{j=1}^n b_j f_0^{(-3^j)}
& = \sum_{j=1}^n b_j f_{3^j}\\
\noalign{\vskip6pt}
& = \sum_{j=1}^n b_j (a_{j-1} r_j -a_{j+1} r_{j+1})\\
\noalign{\vskip6pt}
& = \sum_{j=1}^n (b_j a_{j-1} - b_{j-1} a_j) r_j - b_n a_{n+1} r_{n+1}\ .
\end{split}
\end{equation*}
We deduce that
\begin{equation*}
\begin{split}
&\Big\| \sum_{j=1}^n a_j f_0^{(-3^j)}\Big\|_p
=\|r_1 - a_n a_{n+1} r_{n+1}\|_p \to 1\
\text{ if }\ n\to \infty\ ,\ \text{ and}\\
&\Big\| \sum_{j=1}^n (-1)^{j+1} a_j f_0^{(-3^j)}\Big\|_p
= \Big\|  r_1 + \sum_{i=2}^n (-1)^{i+1} 2a_{i-1} a_i r_i \pm a_na_{n+1}
r_{n+1}\Big\|_p \\
&\hskip1.5truein \sim \bigg( \sum_{i=1}^n a_i^4\bigg)^{1/2}
= \bigg( \sum_{i=1}^n \frac1i\bigg)^{1/2}\ .
\end{split}
\end{equation*}
We can now apply the same arguments as in the proof of Theorem~\ref{thm2.13}
and obtain
$$\Big\| \sum b_j f^{(-3^j)}\Big\|_p \sim
\Big\| \sum b_j f_0^{(-3^j)}\Big\|_p \vee \bigg( \sum |b_j|^p\bigg)^{1/p}\ .
$$
{From} this expression it follows that $(f^{(-3^j)})_{j=1}^\infty$
is basic.

Indeed
$$\Big\| \sum_{j=1}^n b_j f^{(-3^j)}\Big\|_p
\sim\bigg(\sum_{j=1}^n (b_j a_{j-1}-b_{j-1}a_j)^2 + (b_n a_{n+1})^2\bigg)^{1/2}
\vee \bigg( \sum_{j=1}^n |b_j|^p\bigg)^{1/p}\ .$$
Let the right hand expression be equal to 1 with
$$\sum_{j=1}^n (b_j a_{j-1} -b_{j-1} a_j)^2
+ (b_n a_{n+1})^2 =1 .$$
Then if $(b_n a_{n+1})^2 \le 1/2$, for any extension $(b_i)_{i=1}^m$, $m>n$,
the right hand expression is at least $1/\sqrt2$.
If $(b_n a_{n+1})^2 \ge 1/2$ then $b_n > 1/2^{1/4}$ and so
$(\sum_{j=1}^m |b_j|^p)^{1/p} \ge 2^{-1/4}$.

Finally $(f^{(-3^j)})_{j=1}^\infty$ is not unconditional since
$$\Big\| \sum_{j=1}^n a_j f^{(-3^j)}\Big\|_p\sim (\log n)^{1/p}$$
while
$$\Big\|\sum_{j=1}^n (-1)^{j+1} a_j f^{(-3^j)}\Big\|_p\sim (\log n)^{1/2}\ .$$
\qed
\end{example}

The translation problem can, of course, be considered in other
rearrangement invariant function spaces on $\real$.
We end this section with a simple result
in the space $L_p(\real)\cap L_2(\real)$ for $2<p<\infty$.
The norm is given by $\|f\| = \|f\|_p \vee \|f\|_2$ and the space is isomorphic
to $L_p (\real)$ (see e.g., \cite{JMST} for more on this space).

%
%
%

\begin{prop}\label{prop:2.17bis}
Let $2<p<\infty$ and let $(f_i)_{i=1}^\infty$ be an unconditional basis
for $X\subseteq L_p(\real) \cap L_2(\real)$ consisting of translations
of some fixed $f\in L_p(\real)\cap L_2(\real)$.
Then $(f_i)_{i=1}^\infty$ is equivalent to the unit vector basis of $\ell_2$.
\end{prop}

\begin{proof}
As before, by first carefully approximating in both $\|\cdot\|_p$
and $\|\cdot\|_2$ each $f_i$ by a simple dyadic function
$\tilde f_i$ and then choosing a block basis $(g_i)_{i=1}^\infty$ of
$(h_{(i,n)})$ with $|g_i| = |\tilde f_i|$ for all $i$, we obtain:
$(g_i)_{i=1}^\infty$ is equivalent to $(f_i)_{i=1}^\infty$ in
$L_p(\real) \cap L_2(\real)$.

Now $(g_i)_{i=1}^\infty$ is unconditional and semi-normalized in
$L_p(\real) \cap L_2(\real)$ which is isomorphic to $L_p$.
Hence by (1.2), $(g_i)$ admits an upper $\ell_2$-estimate.
Furthermore $(g_i)_{i=1}^\infty$ is unconditional and semi-normalized
in $L_2(\real)$ and thus also admits a lower $\ell_2$-estimate in
$\|\cdot\|_2$ and so in $L_p (\real)\cap L_2(\real)$.
\end{proof}

\section{Discrete versions of the problem}

It remains open if $L_p(\real)$, $1<p<\infty$, admits a basis of
translations of some fixed $f\in L_p(\real)$ (see section~4 for more open
problems).
The examples in section~3 were all integer translations and this leads to
a natural

\begin{quest}\label{quest3.1}
Let $1<p<\infty$.
Is there a set $\Lambda = \{\lambda_n :n\in\nat\}\subseteq \zed$ and an
$f\in L_p(\real)$ so that that $(f_{(\lambda_n)} :n\in\nat)$ is a basis
for $L_p(\real)$?
\end{quest}

The answer is no for $1< p\le 2$ (and of course for $p=1$ by
Theorem~\ref{thm:folk}) by Theorem~\ref{thm:AO-b}.
We also deduce this as a Corollary to Proposition~\ref{prop3.6} below.
The answer is also no for $\Lambda=\zed$ by \cite{AO}.

\begin{prop}\label{prop3.2}
\cite{AO}.
Let $1< p<\infty$.
There is no $\lambda >0$ and $f\in L_p(\real)$ so that $\{f_{(\lambda n)}:
n\in\zed\}$ can be ordered to be a basis for $L_p(\real)$.
\end{prop}

\begin{proof}
%
We will prove a more general result below (see Proposition~\ref{prop3.6} and
Corollary~\ref{cor3.10}).
\end{proof}

We can do a bit better in $L_1$ for certain spaces
$X_1 (f,(\lambda n)_{n\in \zed})$.
By Theorem~\ref{thm:Wi}, $X_1(f,\real) = L_1(\real)$ forces
$\hat f(t)\ne0$ for all $t$.

\begin{lem}\label{lem:minimal}      
Let $f\in L_1(\R)$ with $\widehat f(t)\neq 0$
for all $t$, and let $\Lambda=\{\lambda_n:n\in \nat\}$ be uniformly
discrete. Then $\{f_{(\lambda_n)}\}_{n\in\nat}$ is a non-fundamental
minimal system in $L_1(\R)$.
\end{lem}

\begin{proof}  We use the fact that {\it for a uniformly
discrete $\Lambda$, there exists $a>0$ so that $(e^{i\lambda_n
t})_{n\in\nat}$ is not complete in $C[-a, a]$.}  As pointed out to
us by J. Bruna, this follows from Paley-Wiener theory by
constructing, from an entire function of finite exponential type, a
Paley-Wiener function which vanishes on $\Lambda$. Alternately, this
can be also quickly deduced from Beurling-Malliavin radius of
completeness formula (cf. \cite[section IX D]{Ko}) and the fact that
the uniformly discrete sequences have finite Beurling-Malliavin
density. For convenience of the reader, we present a proof. We
recall the definition of Beurling-Malliavin density $D_{BM}$. For
$\Lambda\subset (0,\infty)$ and $D>0$, a family of disjoint
intervals $(a_k, b_k)$, $0<a_1<b_1<\ldots<a_k<b_k<\ldots
\nearrow\infty$ is called substantial for $D$ if
$$\frac{n_{\Lambda}(a_k, b_k)}{b_k-a_k}>D,\ k=1,2,\ldots,\ \ \ \sum_k
\big(\frac{b_k-a_k}{b_k}\big)^2=\infty,$$ where $n_{\Lambda}(a_k,
b_k)$ is the number points of $\Lambda$ in the interval $(a_k,
b_k)$. Then the density is defined by
$$D_{BM}(\Lambda)=\sup\{D>0:\textrm{there\ exists\ a\ substantial\ family\ for\
D}\}.$$  For a general $\Lambda$, put $D_{BM}(\Lambda)=\max
\{D_{BM}(\Lambda^+), D_{BM}(\Lambda^-)\}$ where
$\Lambda^+=\Lambda\cap \R^+$, $\Lambda^-=(-\Lambda)\cap \R^+$.
Beurling-Malliavin radius of completeness theorem asserts that
$\{e^{i\lambda_n t}:\lambda_n\in\Lambda\}$ is complete in $C[-a,a]$
if and only if $\pi D_{BM}(\Lambda)\ge a$.

Now suppose that $\Lambda$ is uniformly discrete and let
$\delta=\inf \{|\lambda-\lambda'|:\lambda, \lambda'\in\Lambda,
\lambda\neq\lambda'\}>0$. Since
$n_{\Lambda}(a_k,b_k)/(b_k-a_k)<2/\delta$ for all $b_k>a_k>0$, no
$D>2/\delta$ can be substantial for $\Lambda$, and therefore
$D_{BM}(\Lambda)<2/\delta$. Thus, by Beurling-Malliavin theorem,
$(e^{i\lambda_n t})_{n\in\nat}$ is not complete in $C[-b, b]$ for
$b>2/\delta$.

To see the minimality of $\{f_{(\lambda_n)}\}$, suppose to the
contrary that for some $n_0$, $f_{\lambda_{n_0}}\in
[(f_{\lambda_n})_{n\neq n_0}]$ in $L_1(\R)$. Then $\hat
f_{\lambda_{n_0}}(t)= \hat f(t) e^{-i\lambda_{n_0}t} \in [(\hat f(t)
e^{-i\lambda_n t})_{n\neq n_0}]$ in $C_0(\R)$. Now $\hat f(t)\neq 0$
for all $t$, so $e^{-i\lambda_{n_0}t}\in [\{e^{-i\lambda_nt}:n\neq
n_0\}]\subset C[-b, b]$ for all $b>0$. Thus $(e^{-i\lambda_n
t})_{n\neq n_0}$ is complete in $C[-b, b]$ (cf. \cite[Theorem
8, p. 129]{Yo}). This contradicts the fact when $b>a$. Similarly, observe
that $\{f_{(\lambda_n)}\}$ cannot be fundamental in $L_1(\R)$,
indeed otherwise $(e^{-i\lambda_n t})_{n\in\nat}$ would be complete in
$C[-b, b]$ for all $b>0$.
\end{proof}

Note that the assumption $\hat f(t)\ne 0$ for all $t$ is not frivolous
due to Remark~\ref{rem2.7}b).

\begin{prop}\label{prop3.4bis}
Let $f\in L_1(\real)$ with $\hat f(t)\ne0$ for all $t$ and let
$\lambda >0$.
Then $X_1 (f,(\lambda n)_{n\in\zed})$ embeds almost isometrically
into $\ell_1$.
\end{prop}

\begin{proof}
By Corollary~\ref{cor2.5}    
it suffices to show that $(f_{(\lambda n)})_{n\in\zed}$
is a bounded minimal system.
By Lemma~\ref{lem:minimal} it is a minimal system.
Let $g(f)=1$, $g(f_{(\lambda n)})=0$ for $n\in\zed \setminus \{0\}$,
$g\in L_\infty (\real)$.
Then $(f_{(\lambda n)} ,g_{(\lambda n)})_{n\in\zed}$
is a bounded minimal system.
\end{proof}

Proposition~\ref{prop3.2}  generalizes to $\ell_p$-sums of a separable infinite
dimensional Banach space $X$.
Define $\ell_p (X) = \ell_p (\zed,X) = (\bigoplus_{n\in\zed} X)_{\ell_p}$.
For $F= (f_n :n\in\zed) \in \ell_p (X)$ and $k\in\zed$, let
$F^{(k)}$ be $F$ shifted right $k$ times.
Precisely, $F^{(k)} = (f_{n-k})_{n\in\zed}$.

\begin{prop}\label{prop3.3}
Let $X$ be a separable infinite dimensional Banach space, $1\le p<\infty$.
There does not exist $F\in \ell_p(\zed,X)$ so that $\{F^{(k)} :k\in\zed\}$
is a basis for $\ell_p (\zed,X)$ in some order.
\end{prop}

\begin{proof}
Let $1/p + 1/q =1$ and assume for some $F$ that $(F^{(n_i)})_{i=1}^\infty$
is a basis for $\ell_p (\zed,X)$ where $(n_i)_{i=1}^\infty$ is a reordering
of $\zed$.
Let $(G_i)_{i\in \nat}\subseteq \ell_q (\zed,X^*)$ be the biorthogonal
functionals to $(F^{(n_i)})_{i=1}^\infty$.
Choose $i_0$ with $n_{i_0} =0$ and set $G_{i_0} = G = (g_n)_{n\in\zed}$,
with $g_n \in X^*$ for $n\in \nat$.

For $n,m\in\zed$,
$$\langle F^{(n)},G^{(m)}\rangle
= \sum_{k\in\zed} \langle f_{k-n},g_{k-m}\rangle
= \sum_{k\in\zed} \langle f_{k+m-n} ,g_k\rangle
= \langle F^{(n-m)}, G_{i_0}\rangle
= \delta_{(m,n)}\ .$$

Again, from the uniqueness of the biorthogonal functionals to a basis
(for $\ell_p (\zed,X)$), we see that $G_i = G^{(n_i)}$ for all $i\in\nat$.

Choose $j\in \nat$ with
$$\bigg( \sum_{i=j+1}^\infty  \|f_{-n_i}\|^p \bigg)^{1/p}
\le \frac1{2\|G\|}\ .$$

Since $X$ is infinite dimensional, there exists $x\in S_X$ with
$g_{-n_i} (x) =0$ for all $i\le j$.
Set $H = (\delta_{(0,n)} x :n\in\zed) \in \ell_p (\zed,X)$.
Then
$$H = \sum_{i=1}^\infty \langle H,G^{(n_i)}\rangle F^{(n_i)}
= \sum_{i=j+1}^\infty \langle H,G^{(n_i)}\rangle F^{(n_i)}\ .$$
Hence,
\begin{equation*}
\begin{split}
&1=\|x\| =\|H\| =\Big\|\sum_{i=j+1}\langle H,G^{(n_i)}\rangle F^{(n_i)}\Big\|
= \sum_{i=j+1}^\infty \langle g_{-n_i} , f_{-n_i}\rangle \\
\noalign{\vskip6pt}
&\qquad \le \sum_{i=j+1}^\infty \|g_{-n_i}\| \, \|f_{-n_i}\|
\le \|G\| \bigg( \sum_{i=j+1}^\infty \|f_{-n_i}\|^p\bigg)^{1/p} \le \frac12\ ,
\end{split}
\end{equation*}
a contradiction.
\end{proof}

\begin{prob}\label{prob3.5}
Let $2<p<\infty$ and let $X$ be a Banach space with $\dim X\ge2$.
Does there exist $F\in \ell_p (\zed,X)$ and $(\lambda_i :i\in\nat)
\subseteq \zed$ so that $(F^{(\lambda_i)})_{i=1}^\infty$ is a basis for
$\ell_p (\zed,X)$?
What if $\dim X=2$ or if $X = \ell_p$?
\end{prob}

We do not ask the question for $p\le2$ because of the following proposition
which generalizes Proposition~\ref{prop3.3} in that case.

\begin{prop}\label{prop3.6}
Let $1\le p\le 2$ and let $X$ be a Banach space with $\dim (X) \ge2$.
Let $F = (f_i : i\in\zed)\in \ell_p (\zed,X)$.
Then $[\{F^{(n)} :n\in\zed\}]\ne \ell_p (\zed,X)$.
\end{prop}


\begin{cor}\label{cor:AO}
\cite{AO}.
Let $1<p\le 2$, $f\in L_p(\real)$ and $\lambda >0$.
Then $[f_{(\lambda n)} :n\in\zed\}]$ is a proper subspace of $L_p(\real)$.
In particular, no subsequence of $\{f_{(\lambda n)} :n\in\zed\}$ can be
ordered to form a basis for $L_p(\real)$.
\end{cor}

\begin{proof}
We let $F$ denote the Fourier transform on $L_1(\real) + L_2 (\real)$
into the space of measurable functions on $\real$.
$F$ is a bounded linear operator, restricted to $L_1(\real)$ (into
$C_0(\real)$) and when restricted to $L_2(\real)$ (into $L_2(\real)$).
By the Riesz-Thorin interpolation theorem, $F$ is also bounded as a linear
operator from $L_p(\real)$ into $L_q(\real)$ ($1/p + 1/q =1$).
Now since $L_1 (\real) \cap L_2 (\real) \subseteq L_p (\real)$,
$F(L_p(\real))$ is dense in $L_q (\real)$.
For $f\in L_p (\real)$ and $s\in \real$ we have $F(f_s) = e^{-is(\cdot)}F(f)$.
Indeed for $f\in L_1 (\real)$ and $t\in \real$,
\begin{equation*}
\begin{split}
F(f_s)(t) &= \frac1{\sqrt{2\pi}} \int_{-\infty}^\infty e^{ixt} f(x-s)\,dx \\
\noalign{\vskip6pt}
& = \frac1{\sqrt{2\pi}} \int_{-\infty}^\infty e^{-i (u+s)t} f(u)\,dx
= e^{-is} F(f)(t)\ .
\end{split}
\end{equation*}
For a general $f\in L_p (\real)$ the result follows by the standard
density argument.

Let $f\in L_p(\real)$ and $\lambda \in\real$.
If $[\{f_{(\lambda n)} :n\in\zed\}] = L_p(\real)$ then
$[\{ e^{in\lambda (\cdot)} F(f) :n\in\zed\}] = L_q (\real)$.
This implies that $F(f) \ne0$ a.e. and that $[\{e^{in\lambda (\cdot)}:
n\in\zed\}] = L_q (|F(f)|^q\,dx)$.
This in turn implies that all elements $g$ of $L_p (|F(f)|^q\,dx)$ are
$\lambda$-periodic ($g(x) - g(x+\lambda) =0$ a.e.), a contradiction.
\end{proof}

\begin{remark}\label{rem3.7}
For $2<p<\infty$ it is shown in \cite{AO} (Theorem~\ref{thm:AO} above)
that there exists $f\in L_p(\real)$ so that $[(f^{(n)})_{n\in\zed}]
= L_p(\real)$.
\end{remark}

We will use the Fourier transform on the abelian group $\zed$ (see \cite{Ru})
and also assume our spaces to be defined over the complex field.
For $x= (\xi_j) \in \ell_1(\zed)$ we let $\widehat x$ be the function
$$\widehat x :[-\pi,\pi] \to \real\ ,\quad
\widehat x (t) = \sum_{n\in\zed} \xi_n e^{int}\ .$$
It is easy to see that $\widehat x \in C(T)$ when $x\in \ell_1 (\zed)$
(identifying, as usual, the torus $T$ with $[-\pi,\pi]$ by identifying
$\pi$ and $-\pi$).
Also the map
$$\widehat{(\cdot)} : \ell_1 (\zed) \to C(T)\ ,\qquad x\to \hat x$$
is a bounded linear operator of norm 1.
For any $x= (\xi_n : n\in\zed)$
\begin{equation*}
\|\hat x \|_2^2
 = \int_{-\pi}^\pi \Big| \sum_{n\in\zed} \xi_n e^{int} \Big|^2 \,dt
= \int_{-\pi}^\pi \sum_{m,n\in\zed} \bar\xi_m \xi_n e^{i(n-m)t} \,dt
= 2\pi \sum_{n\in\zed} |\xi_n|^2\ .
\end{equation*}
Thus $\widehat{(\cdot)}$ extends to an isometry from $\ell_2 (\zed)$ to
$L_2 (T,\frac1{2\pi}\,dx)$.

Again, by Riesz-Thorin interpolation, the Fourier  transform is a bounded
linear operator from $\ell_p(\zed)$ into $L_q (T)$ for $1\le p\le 2$,
$1/p + 1/q = 1$.

Since $\{\widehat x:x\in\ell_1 (\zed)\}$ is dense in $L_2(T)$, it follows that
the image under the Fourier transform of $\ell_p(\zed)$ is dense in $L_q(T)$.

We also need two lemmas before proving Proposition~\ref{prop3.6}.
The first is an easy exercise in real analysis.

\begin{lem}\label{lem3.8}
Let $\nu\ll \mu$ be two $\sigma$-finite measures on the measure space
$(\Omega,\Sigma)$.
Then for $1\le p<\infty$, if $D\subseteq L_p(\nu) \cap L_p(\mu)$ is dense in
$L_p(\mu)$, it is also dense in $L_p(\nu)$.
\end{lem}

\begin{proof}
Let $\rho$ be the Radon-Nikodym density of $\nu$ with respect to $\mu$.
For $n\in\nat$ set
\begin{equation*}
A_n = \Big\{\omega \in \Omega : \frac1n \le \rho (\omega) \le n\Big\}\ .
\end{equation*}
For $n\in\nat$, it follows that $L_p (\mu |_{A_n}) = L_p (\nu |_{A_n})$.
Also by canonically identifying $L_p (\nu |_{A_n})$ with a subspace
of $L_p(\nu)$, $\bigcup_{n\in\nat} L_p (\nu |_{A_n})$ is dense in
$L_p(\nu)$ and this yields the results.
\end{proof}

\begin{lem}\label{lem3.9}
Let $1\le p\le 2$ and let $x= (\xi_n :n\in\zed)\in \ell_p(\zed)$.
Then  $[(x^{(2n)})_{n\in\zed} ] \ne \ell_p (\zed)$.
\end{lem}

\begin{proof}
Recall  $x^{(n)} = (\xi_{j-n} :j\in\zed)$, for $n\in\zed$.
For $n\in\nat$, $t\in T$ and $z= (\zeta_j :j\in\zed)\in \ell_1 (\zed)$
we have
$$\widehat{z^{(n)}} (t) = \sum_{j\in\zed} \zeta_{j-n} e^{ijt}
= \sum_{\ell\in \zed} \zeta_\ell e^{i(\ell+n)t}
= e^{int} \widehat z\ .$$
By a density argument, we see that for any $x\in \ell_p(\zed)$ and $n\in\zed$,
$\widehat{x^{(n)}} = e^{in(\cdot)} \widehat x$.

Assume, to the contrary, that $[(x^{(2n)})_{n\in\zed}] = \ell_p(\zed)$.
It then follows that $[\{e^{i2n(\cdot)} \widehat x : n\in\zed\}] = L_q(T)$
and thus $\widehat x \ne0$ a.e.
Also that
$$[\{ e^{i2n(\cdot)} : n\in\zed \}]= L_q (T,|\widehat x|^q \,dt)\ .$$
By Lemma~\ref{lem3.8}, this implies that
$$[\{ e^{i2n(\cdot)} :n\in\zed\} ] = L_q (T)\ .$$
Since for any $n\in\zed$,
\begin{equation*}
\begin{split}
& 2\pi \langle e^{i2n(\cdot)} , \chi_{[-\pi,0]} - \chi_{[0,\pi]}\rangle
= \int_{-\pi}^0 e^{i2nt} dt - \int_0^\pi e^{i2nt} dt \\
\noalign{\vskip6pt}
&\qquad = \int_0^\pi (e^{-i2nt} - e^{i2nt})\,dt
= \begin{cases} 0\ ,&\text{if $n=0$}\\
\noalign{\vskip6pt}
\ds -2 \int_0^\pi \sin (2nt) dt =0\ ,&\text{if $n\ne 0$,}
\end{cases}
\end{split}
\end{equation*}
this cannot be true.
\end{proof}

\begin{proof}[Proof of Proposition~\ref{prop3.6}]
After projecting $X$ onto $\ell_p^2$ we see that we may assume $X=\ell_p^2$.
Let $I$ be the obvious isometry between $\ell_p (\zed,X)$ and $\ell_p(\zed)$
denoted
$$(x_j)_{j\in\zed} \longmapsto (y_j)_{j\in\zed}$$
where if $x_j  = (x_{(j,1)}, x_{(j,2)})\in \ell_p^2$ then
$$y_{2j} = x_{(j,1)}, \qquad y_{2j+1} = x_{(j,2)}\ .$$
Then for $(x_j)_{j\in\zed} \in \ell_p (\ell_p^2)$,
$(x^{(n)})_{n\in\zed} = (I(x)^{2n})_{n\in\zed}$ and the result follows
from Lemma~\ref{lem3.9}.
\end{proof}

\begin{rem}
As noted above by
the results of \cite{AO} in section~4 we cannot hope to prove that given
$f\in L_p (\real)$, $2<p<\infty$, $[(f^{(n)})_{n\in\zed}]\ne L_p(\real)$.
Nevertheless, by dualizing Proposition~\ref{prop3.6}, we
have the following
\end{rem}

\begin{cor}\label{cor3.10}
Let $X$ be a Banach space with $\dim (X) \ge2$ and let $2\le p<\infty$.
Let $F= (f_i)_{i\in\zed} \in \ell_p (\zed,X)$.
Then $\{ F^{(n)} :n\in\zed\}$ is not a basis for $\ell_p(\zed,X)$
under any ordering.
\end{cor}

\begin{proof}
Assume that $F\in (f_i)_{i\in\zed} \in \ell_p (\zed,X)$ and that
$(n_s)_{s\in\nat}$ is an ordering of $\zed$ so that $(F^{(n_s)})_{s=1}^\infty$
is a basis for $\ell_p (\zed,X)$.
Let $(G_s)_{s=1}^\infty\subseteq \ell_q (\zed,X^*)$ be the biorthogonal
functionals of $(F^{(n_s)})_{s=1}^\infty$.
Set $G = (g_j)_{j\in\zed} = G_1$.
We let $G^{(m)} = (g_{j-m})_{j\in\zed}$, as usual.
For $s,t\in \nat$ and $m\in\zed$ we have
\begin{equation*}
\begin{split}
\langle F^{(n_s)}, G^{(n_t)}\rangle
& = \sum_{j\in\zed}\langle f_{j-n_s}, g_{j-n_t}\rangle
= \sum_{k\in \zed} \langle f_{k+n_t - n_s} , g_k\rangle\\
\noalign{\vskip6pt}
& = \langle F^{(n_s - n_t)},G_1\rangle
= \begin{cases} 1\ ,&\text{if $n_s - n_t =n_1$}\\
\noalign{\vskip6pt}
0\ ,&\text{if $n_s- n_t \ne n_1$.}
\end{cases}
\end{split}
\end{equation*}
As before, we see that $G_s = G^{(n_s-n_1)}$.
In particular, span$\{ G^{(n)} :n\in\zed\}$ is $w^*$-dense in $\ell_q(X^*)$.
Let $E$ be a two dimensional subspace of $X$ and let $P$ be a projection
of $X$ onto $E$.
Let $Q : \ell_p (\zed,X)\to \ell_p (\zed,E)$ be the projection given
by $Q(H) = (P(h_i))_{i\in\zed}$.
It follows that span$(G^{(n)} |_{\ell_p(\zed,E)} )_{n\in\zed}$ is
$w^*$ dense in $\ell_q (\zed,E^*)$ and hence norm dense (the latter is
reflexive).
This contradicts Proposition~\ref{prop3.6}.
\end{proof}

\section{Results from the literature and open problems}

We first cite some more known results from the literature.

\begin{thm}\label{thm4.5} \cite[Theorem 5.1(b)]{DH}.
Let $g^{(1)}, g^{(2)}, \ldots, g^{(m)} \in L_2(\real^d)\cap L_1(\real^d)$
and let $\Gamma_1, \Gamma_2,\ldots,\Gamma_m \subset \real^d$ be countable.
Then $\{ g_{(\lambda)}^{(i)} : i = 1,2,\ldots,m,\, \lambda \in \Gamma_i\}$
cannot be ordered to be a Schauder basis of $L_2(\real^d)$.
\end{thm}

\begin{thm}\label{thm4.6}
{\rm (\cite{ER} and \cite{Ro}, cf.\ \cite[Theorem 9.18]{H})}
If $g\in L_p(\real^d)$, $g\ne 0$, and $1\le p\le  \frac{2d}{d-1}$
then the functions $\{g((\cdot) - \alpha_k) : k=1,2,\ldots,N\}$ are
linearly independent for any $N\in\nat$ and any collection
$(\alpha_k)_{k=1}^N \subseteq \real^d$ of distinct points.

If $\frac{2d}{d-1} <p\le\infty$, then for $N\in\nat$ there exists
$0\ne g\in L_p (\real^d)$ and distinct points $(\alpha_k)_{k=1}^N\subseteq
\real^d$ so that $\{ g((\cdot)-\alpha_k) : k=1,2,\ldots,N\}$ is linearly
dependent.
\end{thm}

Our last cited result requires some notation.
For $\Lambda \subseteq \real$ let
$\E (\Lambda) = \text{span}\{ e^{i\lambda (\cdot)} :\lambda\in\Lambda\}$.
Let $R(\Lambda) = \sup \{\rho>0 :\E(\Lambda)$ is dense in
$C[-\rho,\rho]\}$.
Recall, $\Lambda \subseteq \real$ is {\em discrete} if it has no
accumulation points.

\begin{thm}\label{thm4.7} \cite[Theorem 1]{BOU}.
Let $\Lambda \subseteq \real$ be discrete.
There exists $f\in L_1 (\real)$ so that $[\{f_{(\lambda)} :\lambda\in\Lambda\}]
= L_1 (\real)$ if and only if $R(\Lambda) = \infty$.
\end{thm}

Finally we list some problems that remain open.
The main one is

\begin{prob}\label{prob4.8}
Let $1<p<\infty$.
Does there exist $f\in L_p(\real)$ and $\Lambda\subseteq  \real$ so that
$\{f_{(\lambda)} : \lambda\in\Lambda\}$ can be ordered to be a basis
for $L_p (\real)$?
Can we find $f$ and a uniformly discrete set $\Lambda$
so that $\{f_{(\lambda)}:\lambda \in\Lambda\}$
can be ordered to be a frame for $L_p(\real)$?
\end{prob}

By identifying $L_p(\real)$ with $L_p[0,1]$ we have a more general
version of the basis problem in \ref{prob4.8}.

\begin{prob}\label{prob4.9}
Does there exist a normalized basis $(f_n)_{n=1}^\infty$ for $L_p[0,1]$,
$1<p<\infty$, so that for all $0<b<1$,
$$\sum_{n=1}^\infty \| f_n |_{[0,b]}\|^p < \infty\ ?$$
If $4<p < \infty$ can we find such $f_n$'s which form an unconditional
basis for $L_p [0,1]$?
\end{prob}

\begin{prob}\label{prob4.10}
Let $4 <p<\infty$.
Does there exist $f\in L_p(\real)$ and $\Lambda \subseteq \real$ so that
$(f_{(\lambda)})_{\lambda\in\Lambda}$ is an unconditional basis for
$L_p(\real)$?
\end{prob}


We can also raise questions asking for less and here is one such question.

\begin{prob}\label{prob4.11}
Let $1< p<4$.
Does there exist $f\in L_p(\real)$ and a uniformly
discrete set $\Lambda\subseteq \real$
so that $[\{f_{(\lambda)} : \lambda\in \Lambda\}] \subseteq L_p(\real)$
contains an isomorph of $\ell_2$ and $(f_{(\lambda)})_{\lambda\in\Lambda}$
can be ordered to be a basic sequence (or a frame)?
\end{prob}

\begin{prob}\label{prob4.12}
Let $\Lambda \subseteq \real$ be uniformly discrete and $f\in L_1(\real)$.
Does $X_1 (f,\Lambda)$ embed into $\ell_1$?
\end{prob}



\begin{thebibliography}{CDOSZ}
\frenchspacing

\bibitem[AOd]{AOd}
D. Alspach and E. Odell,
{\em $L_p$ spaces},
in ``Handbook of Geometry of Banach Spaces'', Vol.1,
W.B.~Johnson and J.~Lindenstrauss, eds.,
Elsevier, Amsterdam (2001), 123--159.


\bibitem[AO]{AO} A.~Atzmon and A.~Olevskii,
{\em Completeness of integer translates in function spaces on $\R$},
J. of Approx. Theory { \bf 87} (1996) 291--327.

\bibitem[BOU]{BOU}  J.~Bruna, A. Olevskii and A. Ulanovskii,
{\em Completeness in $L_1(\R)$ of discrete translates},
Rev. Mat. Iberoamericana  {\bf 22} no. 1,  (2006),  1--16.

\bibitem[CDOSZ]{CDOSZ}
P.G. Casazza, S.J. Dilworth, E. Odell, Th. Schlumprecht and A. Zs\'ak,
{\em Coefficient quantization for frames in Banach spaces},
J. Math. Anal. Appl. {\bf 348} (2008), 66--86.

\bibitem[CHL]{CHL}
P.G. Casazza, D. Han and D.R. Larson,
{\em Frames for Banach spaces, the functional and harmonic analysis
of wavelets and frames},
(San Antonio, TX, 1999),
Contemp. Math. {\bf 247} (1999), 149--182.



\bibitem[DH]{DH} B.~Deng and C. Heil,
{\em Density of Gabor Schauder bases}, in: Wavelet Applications
in Signal and Image Processing VIII (San Diego, CA, 2000), A. Aldroubi,
A. Lane, and M. Unser, eds., Proc. SPIE {\bf 4119}, SPIE, Bellingham, WA,
 (2000),  153--164.

\bibitem[ER]{ER} G.~Edgar and J.~Rosenblatt,
{\em  Difference  equations over locally compact abelian groups},
Trans. Amer. Math. Soc. {\bf  253} (1979),  273--289.


\bibitem[H]{H} C.~Heil,
{\em  Linear independence of finite Gabor systems}
in ``Harmonic Analysis and Applications,''
A volume in honor of John J. Benedetto,
Birkhauser, Boston, (2006),  171--206.

\bibitem[J]{J}
W.B.~Johnson,
{\em On quotients of $L_p$ which are quotients of $\ell_p$},
Compositio Math. {\bf 34} (1) (1977), 69--89.

\bibitem[JL]{JL}
W.B.~Johnson and  J. Lindenstrauss,
{\em Basic concepts in the geometry of Banach spaces},
in ``Handbook of Geometry of Banach Spaces'', Vol.1,
W.B.~Johnson and J.~Lindenstrauss, eds.,
Elsevier, Amsterdam (2001), 1--84.

\bibitem[JO]{JO}
W.B.~Johnson and E.~Odell, {\em Subspaces of $L_p$ which embed into $\ell_p$},
Compos.  Math. {\bf 28} (1974), 37 --49.


\bibitem[JMST]{JMST}
W.B.~Johnson, B. Maurey, G. Schechtman and L. Tzafriri,
{\em Symmetric structures in Banach spaces},
Mem. Amer. Math. Soc. {\bf 19}(217) (1979) v+298.

\bibitem[KP]{KP}
M.I. Kadets and A. Pe{\l}czy\'nski,
{\em Bases, lacunary sequences and complemented subspaces in the spaces $L_p$},
Studia Math. {\bf 21} (1961/62), 161--176.

\bibitem[KW]{KW}
N.J. Kalton and D. Werner,
{\em Property $(M)$, $M$-ideals, and almost isometric structure of Banach
spaces},
J. Reine Angew. Math. {\bf 461} (1995), 137--178.

\bibitem[Ko]{Ko}
P. Koosis,
{\em The logarithmic integral II}, Cambridge
Studies in Advanced Mathematics, 21. Cambridge University Press,
Cambridge, 1992.



\bibitem[LT]{LT}
J. Lindenstrauss and L. Tzafriri,
{\em Classical Banach Spaces I: Sequence Spaces},
Ergebnisse der Mathematik und ihrer Grenzgebiete 92,
Springer-Verlag, Berlin (1977).

\bibitem[Ol]{Ol}
A.~Olevskii,
{\em Completeness in $L_2(\R)$ of almost integer translates},
C.~R.~Acad.~Sci.~Paris {\bf 324 } (1997), 987--991.


\bibitem[OZ]{OZ}
T. E. Olson and R. A. Zalik,
{\em  Nonexistence of a Riesz basis of translates, in: Approximation Theory},
Lecture Notes in Pure and Applied Math.,  {\bf 138},
Dekker, New York (1992) 401--408.

\bibitem[R]{R}
J. Rosenblatt, private communication.

\bibitem[Ro]{Ro}
J.~ Rosenblatt,
{\em Linear independence of translations},
J. Austral. Math. Soc. (Series A) {\bf 59} (1995), 131-133.

\bibitem[Ru]{Ru} W.~Rudin,
{\em Fourier Analysis on Groups},
Interscience Publisher, John Wiley \& sons (1967).

\bibitem[S]{S}
G. Schechtman,
{\em A remark on unconditional basic sequences in $L_p$ $(1<p<\infty)$},
Israel J. Math. {\bf19} (1974), 220--224.

\bibitem[Wi]{Wi} N.~Wiener,
{\em The Fourier integral and certain  of its applications},
Cambridge University Press,
 Cambridge, 1933 reprint: Dover, New York, 1958.

\bibitem[Yo]{Yo} R. M. Young,
{\it An introduction to nonharmonic Fourier series},
Pure and Applied Mathematics, 93. Academic Press, Inc., New
York-London, 1980.



\end{thebibliography}
\end{document}